\documentstyle[12pt]{article}
\newtheorem{df}{Definition}
\newtheorem{teo}{Theorem}
\newtheorem{prop}{Proposition}
\newtheorem{lem}{Lemma}
\newtheorem{cor}{Corollary}
\newtheorem{rem}{Remark}
\newtheorem{nota}{Notation}
\newtheorem{ex}{Example}
\title{Focal Loci of Algebraic Varieties I}
\author{Fabrizio Catanese - Cecilia Trifogli}
\date{april 2000}

\begin{document}

\maketitle

\begin{abstract}
The focal locus $\Sigma_X$ of an affine variety $X$ is  roughly
speaking the
(projective)  closure of the
set of points $O$ for which there is a smooth point  $x \in X$
and a circle with centre $O$ passing through $x$
which osculates $X$ in $x$.
Algebraic geometry interprets the focal locus as the branching locus
of the endpoint map $\epsilon$ between the Euclidean normal bundle $N_X$
and the projective ambient space ($\epsilon$ sends the normal vector $O-x$
to its endpoint $O$),
and in this paper we address two general problems :

1) Characterize the "degenerate"
case where the focal locus is not a hypersurface

2) Calculate, in the case where   $\Sigma_X$ is a hypersurface,
its degree (with multiplicity)

\end{abstract}

\section{Introduction}

The goal of the present paper is to introduce a general
theory of focal loci of algebraic varieties in Euclidean space.

The theory of focal loci was classically considered only for plane curves
and surfaces in 3-space ( cf. [Coolidge] , [Salmon-Fiedler]),
and Hilbert himself lectured in the Winter Semester 1893-94 at
the University of G\"ottingen on the focal loci of curves and surfaces
of  degree two in 3-space.

Recently  the theory was considered in ([Fantechi], [Trifogli]) for
the respective cases of plane curves and hypersurfaces.

We would like to first briefly present the relevant concepts.

Usually the focal locus of a submanifold $X$ ( cf. [Milnor],
¤ 6, pp. 32-38, or also [D-F-N], vol. II ¤11, sections 2-3)
is defined in
Euclidean differential geometry as either the locus of centres of
principal curvatures, or, more geometrically, as the locus
where the infinitely near normal spaces intersect each other.
Equivalently, the focal locus can also be defined as
 the complement of the set of points
$p$ such that the square of the distance function from $p$
induces a local Morse function on $X$, or also as the union of the
singular points of the parallel varieties to $X$.

To make the definition algebraic, one picks up the second
geometrical definition, where  the
 notion of length is not needed,
 just the notion of orthogonality is sufficient.

To explain this in more detail, let us consider (complex)
affine space as the complement of
a hyperplane ( the "hyperplane at infinity") in projective space.
In the hyperplane at infinity ${\bf P_{\infty}}$,
we give a non degenerate quadric $ Q_{\infty}$.

These data allow , for each projective linear subspace $L$ ,
 to define the orthogonal
 $L^{\perp_{x}}$ to $L$ in a point $x$ as the join of $x$
with the "orthogonal
direction" to $L$ ( this is the subspace of ${\bf P_{\infty}}$
given by the polar of $L \bigcap {\bf P_{\infty}}$ with respect
to  $ Q_{\infty}$).

Given now an irreducible algebraic variety
$X^{n}_{d} \subset {\bf P^{m}}$,
of dimension $n$ and degree $d$ and not contained in the hyperplane
at infinity ,
 for each smooth point $x \in X -{\bf P_{\infty}}$
we define the normal space $N_{x}(X)$
as the orthogonal in $x$ of the projective
tangent space to $X$ at $x$. The condition that $x$ is a point in affine space
ensures that $N_{x}(X)$  has the correct dimension $m-n$.

The normal variety  $N_{X}$ is then defined as the irreducible
algebraic set  in
 ${\bf P}^{m}\times {\bf P}^{m}$, closure
of the  set $N^{good}_{X}$ consisting of the pairs $ (x,y)$
where $x$ is
a smooth point of $X$,
 $x \in X -{\bf P_{\infty}}$ and $y \in N_{x}(X)$.

Clearly, $N_{X}$ is a projective variety of dimension $m$
and the second projection induces a map $\pi$ whose image
is the closure of the union of the normal spaces to the
smooth points of  $ X -{\bf P_{\infty}}$. Observe moreover that
$N^{good}_{X}$ is a projective bundle over $X-{\bf P_{\infty}}-
Sing(X) $, in particular $N^{good}_{X}$ is smooth of dimension
$m$ : therefore we can consider the ramification locus
$Y^{good}_{X}$ of $ \pi: N^{good}_{X} \rightarrow {\bf P}^{m}$,
and we define the ordinary ramification locus as the closure
$Y_{X}$  of $Y^{good}_{X}$.

Defining the good focal locus
as $\Sigma^{good}_{X} = \pi ( Y^{good}_{X})$, and the focal
 locus $ \Sigma_{X}$ as the closure of $\Sigma^{good}_{X}$ ( thus
$ \Sigma_{X}$ is contained in the branch locus of
$\pi : N_{X} \rightarrow  {\bf P}^{m}$),
we have a priori at least four cases:

\begin{itemize}

\item
1) $\pi : N_{X} \rightarrow  {\bf P}^{m}$ is not dominant :
in this case we say that the variety $X$ is  isotropically focally
 degenerate (for short : isotropically degenerate),
and observe that the focal locus $ \Sigma_{X}$ of $X$ is
then simply the image of $ \pi$ ( whence, $ \Sigma_{X}$ is an irreducible
variety in this case!).

\item
2) $\pi : N_{X} \rightarrow  {\bf P}^{m}$ is dominant ,
but the focal locus $ \Sigma_{X}$ (respectively, the branch locus of
$\pi : N_{X} \rightarrow  {\bf P}^{m}$) has dimension
 at most $m-2$ : in this case we say that $X$ is strongly
focally degenerate ( respectively, completely strongly
focally degenerate).

\item
3) $\pi : N_{X} \rightarrow  {\bf P}^{m}$ is dominant ,
 whence surjective , and
the focal locus $ \Sigma_{X}$ "is not
a hypersurface", in the sense that not every component $Z$
of the ordinary ramification divisor $Y_{X}$ (closure of $Y^{good}_{X}$)
 maps to a hypersurface.
In this case we shall say that $X$ is weakly focally degenerate.
We shall moreover say that we have the vertical case if $Z$
 does not dominate $X$.

\item
4) When none of the above occurs, in particular
$\pi : N_{X} \rightarrow  {\bf P}^{m}$ is
 surjective , and the focal locus $ \Sigma_{X}$ is
a hypersurface, we shall say that $X$ is focally non degenerate.
In this case,  defining the focal hypersurface as a divisor,
consisting as the image of the ramification divisor $Y_{X}$
with multiplicities (if $Y_{X} = \Sigma_{i=1,..k}
 n_i Y_i$ , and $ d_i : = degree (Y_i \rightarrow \pi (Y_i)$,
then, setting $\Sigma_i : = \pi (Y_i)$, we get
$ \Sigma_{X}: =  \Sigma_{i=1,..k} d_i   n_i \Sigma_i$),
the main problem is to describe $\Sigma_X$.

\end{itemize}

The first main result of this paper consists in
calculating the degree (with multiplicity) of the focal hypersurface
under a certain hypothesis upon $X$, which we call of being
"orthogonally general", and which ensures that $X$ is focally
non strongly degenerate if it is not a linear subspace.
This concept is important because, if $X$ is smooth and
not a linear subspace,
then  for a general projectivity
$g$ the translate  $g(X)$ of $X$ by $g$ satisfies this condition
whence it is not focally strongly degenerate and we have a divisor $\Sigma_X$.
 The hypothesis that $X$ be "orthogonally general"
is indeed very easy to verify since it simply amounts to
 three requirements: the smoothness
 of $X$,  plus the two
general position properties that $X$ be transversal
to ${\bf P_{\infty}}$, respectively to  $Q_{\infty}$.

More precisely, we have the following Theorem\\

{\bf Theorem \ref{orthgen}}
{\it Let $X\subset {\bf P}^{m}$ be a variety of dimension $n \geq 1$
which is orthogonally general.  Then
$dim\, \Sigma_{X}<m-1 \Leftrightarrow X$ is a linear space.
If $X$ is a linear space, $\Sigma_X$ is a linear space of
dimension equal to $codim\, X-1$. }\\

One can ask in the above theorem whether one can replace the condition
$dim\, \Sigma_{X}<m-1$ (i.e., that $X$ be strongly focally degenerate)
by the weaker condition that $X$ be focally degenerate.

As a corollary of the full description given in Theorem 3 of the
focally degenerate varieties, it turns out that if $X$ is an orthogonally
general and focally degenerate variety, then either $X$ or
$X_\infty$ should be a developable
variety rather explicitly described, but we have not yet had the time
to look at the existence question for such very special varieties.

It is rather clear ( e.g., from the case of plane curves)
that the condition of being orthogonally general is a sufficient
but not necessary condition in order that $X$ be non focally strongly
degenerate. When $X$ is non orthogonally general, but
focally non degenerate, what happens is that the  degree of the focal
divisor can drop ( in this case, for plane curves we have
Pl\"ucker type formulae, cf. [Fantechi]).

Naturally, what we have said insofar opens a series of problems.
To some of them we give an answer in the present paper, to some
others we hope to return in a sequel to this paper :

\begin{itemize}
\item 1) Try to completely classify the  focally isotropically
degenerate varieties.
In section 7 we give a structure Theorem ( Theorem 4)
stating that the isotropically focally degenerate hypersurfaces
are exactly the isotropically developable hypersurfaces.
We observe thus that there are plenty of intriguing
examples already in the case of surfaces in 3-space: these
are obtained as the tangential developable surface of
any space curve whose tangent direction is always an isotropic
vector.
We give moreover a description in section 8, Theorem 5,
 of the general case,
in terms of the inverse focal construction applied to the
focal variety $\Sigma$ and to an algebraic function $r$ on
 $\Sigma$.  We get thus
an implicit classification of these varieties as
developable varieties, but  for this we need to start with a variety
$\Sigma$ whose normal spaces are totally isotropic,
and the function $r$ must also satisfy a suitable condition.

\item 2) Try to classify the weakly and the strongly
focally degenerate varieties.
In section 6 we give a complete classification for the
weakly focally degenerate varieties, showing in Theorem 3
how they can be divided into some"primitive" classes
( cases 1), 2), 6), 7)) and some "derived" classes (cases 3),4),5)),
related for instance by some tangential conditions to
some primitive focally degenerate varieties.
 The primitive focally degenerate varieties
can be described starting from  fibrations in spheres or in affine spaces
"around" the degenerate
component $\Sigma$ of the focal locus.

The question of classifying the strongly focally degenerate
varieties seems harder.

\item 3) Determine whether for a general projective deformation of
$X$ the focal hypersurface is reduced of degree equal to the virtual
degree, and moreover  answer more specific questions such as :

\item 3a)  can we also obtain that for a general deformation
the focal hypersurface has generic Lagrangian singularities ?

\item 3b) can we obtain the above good properties for the
focal hypersurface $ \Sigma_{gX}$ of a general translate
$g X$ of $X$ by a general projectivity $g$ ?

\end{itemize}

Concerning the first problem, the situation seems to us rather hard
(although quite interesting) as soon as the dimension of the
ambient space grows: for instance, whereas a  focally isotropically
degenerate plane curve $C$ is necessarily a line through a cyclic
point ( these are the two points of $Q_{\infty}$ , satisfying
the equations $ z= x^2 + y^2 = 0$ ), in the case of  a  surface
in 3-space we obtain the tangential developable of a space curve $C$
which is "isotropic" in the following sense:
 $C$ is just a curve such that any of its tangent lines $L$
has the property that $L$ intersects ${\bf P_{\infty}}$
in a point of $Q_{\infty}$. Therefore, if we write  the point
of the curve $C$ as a vector function $ x(t)$ of a parameter $t$,
we just have to solve the differential equation

$Q_{\infty} ( dx/dt) = 0 $.

Thus such a curve $C$ yields a  curve $ \Gamma$
in $Q_{\infty}$ parametrizing the projective tangent lines to $C$,
 and the question reduces to: for which $\Gamma$ can
one find an algebraic integral ? (however, since the ring of
polynomials in $t$ is stable by $d/dt$ , the above observation
easily allows us to construct a lot of  focally isotropically
degenerate surfaces, which are tangential surfaces of rational
space curves, cf. Example 10).

In higher dimension, as we already remarked, Theorem 5
partially reduces
the quest to the search of varieties with totally isotropic
normal spaces.

Turning to the other problems, the situation is clear for
the plane curves (cf. [Fantechi]) : the only focally
degenerate plane curves, which are not lines,
are the circles (conics through the two cyclic points),
and moreover, for an irreducible plane curve $C$ the map of
$C$ to the focal curve $ \Sigma_{C}$ is non birational
exactly for a well classified class of curves (by the way,
Fantechi shows that this class is non empty, contrary to a
statement made in [Coolidge]).

As we said, we characterize (cf. Theorem 3) the weakly focally degenerate
varieties, distinguishing six essentially different cases :

\begin{itemize}
\item
two vertical cases, where the exceptional component $R$ of $Y_X$
does not dominate $X$, but is instead the restriction of the normal
bundle $N_X$ to a divisor $ X"$. In both cases, $ X"$ is focally
degenerate, and the focal degeneracy of $X$ is determined by the
first order neighbourhood of $X$ along $X"$ (see Theorem 3 for
more details).

\item
 the case where
$X$ consists of a family of $(m-1-a)$-dimensional spheres parametrized by
the $a$-dimensional degenerate component $\Sigma$ of the focal locus: this
family is  moving according
to a simple differential equation which can be explicitly solved,
and it turns out that we get  a family of spheres each obtained
as the intersection of the big sphere with centre $ O \in \Sigma$
with an affine subspace orthogonal to the tangent space to
$\Sigma$ in $O$.
\item
The case where $X$ is a "transversal" divisor in a focally isotropically
degenerate variety.
\item
The asymptotic case, i.e., the case where $\Sigma$ lies at infinity,
and then $X$ is a developable variety whose intersection $X_{\infty}$
with the hyperplane at infinity "is" the dual variety of $\Sigma$
in ${\bf P_{\infty}}$. In this case there is another simple process,
 called the "asymptotic inverse focal construction", describing
$X$ in terms of the data of $\Sigma$ and of an algebraic function $r(s)$
on $\Sigma$.
\item
The isotropically asymptotic case, where $\Sigma$ lies at infinity, and
a component $\Delta$ of $X_\infty$ is projectively isotropically
 degenerate. This case
is characterized by the property that $\Delta \subset X_\infty$
be obtained via the
isotropic projective inverse focal construction, starting from
$\Sigma, r(s)$ satisfying suitable conditions.

\end{itemize}

The characterization given in Theorem 3
(where also  the case of the focally
 isotropically degenerate varieties is considered) is expressed
in terms of the "inverse focal construction",
which, starting from a variety $\Sigma$ of dimension $a$, and an algebraic
function $r(s)$ on $\Sigma$, considers the union $X'$ of the
 family of spheres each obtained
as the intersection of the big sphere with centre $ O \in \Sigma$
and radius equal to the square root of $r(s)$
with an affine subspace orthogonal to the tangent space to
$\Sigma$ in $O$, and whose position is determined by the
differential of the function $r(s)$.

It turns out that for the focally isotropically
degenerate varieties the above spheres degenerate to
 affine spaces,and $X$ equals $X'$,
 whereas in the case where these spheres have
the right dimension $m-1-a$  $X'$ is focally degenerate.

For hypersurfaces in  higher dimensions the second author
( [Trifogli]) showed that the focal hypersurface
of a general hypersurface  is reduced (indeed  that this holds
for a general diagonal
hypersurface, i.e., for a translate of
the Fermat hypersurface by a projectivity in the diagonal torus).

Concerning problem 3a), this is a global problem which is however
related to a local problem which has been extensively studied:
the theory of Lagrangian singularities.
In fact the Normal variety $N_X$ is a Lagrangian variety for the
symplectic form on the product ${\bf A^m} \times {\bf A^m}$
which is associated to  $Q_{\infty}$, namely
$^t x Q_{\infty} y -  ^t y Q_{\infty} x  $, and the second projection
is also  Lagrangian ( cf. [Arnold et al.]).

Partial results concerning problem 3a) have  been obtained
by the second author for surfaces in 3- space( [Trif2]).

\section{Notation}

$V'$ :=  a fixed vector space of dimension $m$

$V$ := the vector space $V=V'\bigoplus {\bf C}$

${\bf P}(V)={\bf P}^{m}$ := the projective space whose points
correspond to the 1-dimensional vector subspaces of $V$

 ${\bf P}(V')={\bf P_{\infty}} \subset {\bf P}(V)$
$(\cong {\bf P^{m-1}})$ the complement
of the affine space ${\bf P}(V) -{\bf P_{\infty}} \cong V'$.

$X^{n}_{d} \subset {\bf P^{m}}$ a quasi-projective
algebraic variety of dimension
$n$ and degree $d$ which does not lie at infinity , i.e.,

$X^{n}_{d} \not\subset {\bf P_{\infty}}$

$  Q_{\infty}$ := a non degenerate quadratic form on $V'$,
 yielding an isomorphism

$Q:V'\stackrel{\cong}{\rightarrow} (V')^{\vee}$.

By slight abuse of notation, the corresponding quadric

$  Q_{\infty} \subset{\bf P_{\infty}} $.

$W$ := a vector subspace  of $V'$,

$Ann(W)$ := the vector subspace of $V'$ which is the
orthogonal space of $W$ with respect to the
quadratic form $  Q_{\infty}$

$L'$ :=  ${\bf P}(W) $ a linear subspace at infinity

 ${L'} ^{\perp}$ := ${\bf P}( Ann (W)) $ , the polar subspace
of $L'$ , also called the orthogonal direction to $L'$

$L   \subset {\bf P}(V)$ := a projective linear
subspace ,  $L \bigcap {\bf P_{\infty}}$ the direction of $L$

 ${L} ^{\perp_{x}}$ := the orthogonal to $L$ in $x$, defined
as the smallest linear subspace containing $x$ and the
orthogonal direction of $L$ ( i.e., the polar of
$L \bigcap {\bf P_{\infty}}$) .

\section{"Normal Bundle" in Euclidean Setting}

In this section, we shall consider a smooth quasi- projective variety
$X^{n}_{d} \subset {\bf P^{m}}$ and we shall define its
 projective normal variety  $ N_{X} \subset
{\bf P}^{m}\times {\bf P}^{m}$, and its Euclidean Normal sheaf
${\cal N}_X$.

Under some assumptions that we are going to specify, the first
projection of the normal variety  $ N_{X}$ to $X$
yields a projective bundle over $X$, which is the projectivization
of the Euclidean Normal sheaf :

$ N_{X} = {\bf P}({\cal N}_{X})\subset
{\bf P}(V \bigotimes {\cal O}_{X}) \subset
 {\bf P}(V \bigotimes {\cal O}_{{\bf P}^{m}}) =
{\bf P}^{m}\times {\bf P}^{m}$.

Start from the Euler sequence\\

$(1)$ $0\rightarrow {\cal O}_{\bf P}(-1)\rightarrow V\bigotimes
{\cal O}_{\bf P}\rightarrow T_{\bf P}(-1)\rightarrow 0$ :\\

setting ${\cal L} ={\cal O}_{X}(-1)$, the restriction to $X$ of
>the Euler sequence and the inclusion of the
tangent bundle of $X$ in the restricted tangent bundle
of ${\bf P}^{m}$ define the bundle $ {\widehat T}_{X}(-1) $
whose projectivization is the projective tangent bundle to $X$.

We get thus two  exact sequences, the second included into the first:\\

$(2)$   $\begin{array}{l}
0 \rightarrow {\cal L}  \rightarrow V\bigotimes
{\cal O}_{X} \rightarrow T_{\bf P}(-1)\bigotimes {\cal O}_{X}
 \rightarrow 0\\
0  \rightarrow  {\cal L}  \rightarrow  {\widehat T}_{X}(-1)
 \rightarrow
 T_{X}(-1)  \rightarrow  0
\end{array}$\\

$\underline {Assumption\, 0= smoothness}$: $X$ is smooth,
 whence $T_{X}$ and
${\widehat T}_{X}$ are subbundles.\\

 Recalling that $V=V'\bigoplus {\bf C}$,we state the further

$\underline {Assumption\, 1}$ ( = transversality of the intersection
$X \bigcap {\bf P_{\infty}}$
 with the
hyperplane at infinity) : ${\overline T}_{x}:=V'\cap
{\widehat T}_{x}$ is a hyperplane in ${\widehat T}_{x}$ $\forall
x\in X$.\\

This means that we have two more  exact sequences\\

$(3)$   $\begin{array}{l}
0 \rightarrow V'\bigotimes {\cal O}_{X} \rightarrow V\bigotimes
{\cal O}_{X} \rightarrow {\cal O}_{X}
 \rightarrow 0\\
0  \rightarrow  {\overline T}_{X}(-1)  \rightarrow
{\widehat T}_{X}(-1)  \rightarrow
{\cal O}_{X}  \rightarrow  0
\end{array}$\\

At this stage we can define the bundle of normal directions
${\cal N}'_{X}$ as a twist of the
annihilator of  ${\overline T}_{X}$.

  We define it through the
  exact sequence

$(4)$ $0 \rightarrow {\cal N}'_{X}(-1) \rightarrow
V'\bigotimes {\cal O}_{X} \cong (V')^{\vee} \bigotimes {\cal O}_{X}
\rightarrow ({\overline T}_{X}(-1))^{\vee}\rightarrow 0$\\

In order to obtain a projective normal bundle from
the bundle of normal directions we need a last

$\underline {Assumption\, 2}$ ( = transversality of $X$ with
$Q_{\infty}$) : The natural map

 ${\cal L} \bigoplus {\cal N}'_{X}(-1)  \rightarrow V
\bigotimes {\cal O}_{X}$ is a
bundle embedding, thus its image ${\cal N}_{X}(-1)$ is a subbundle of
$V\bigotimes {\cal O}_{X}$, isomorphic to ${\cal L}
\bigoplus {\cal N}'_{X}(-1)$\\

We notice thus that if assumption 2) holds, then ${\cal N}_{X}
\cong {\cal O}_{X} \bigoplus {\cal N}'_{X}$

\begin{df}
$X$ is said to be orthogonally general if it satisfies Assumptions
$0-2$ above.
\end {df}

\begin{rem}
For every algebraic variety $X\subset {\bf P^{m}}$ which is
not contained in the hyperplane at infinity there is a
maximal nonempty
Zariski open set $U$ of $X$ which is orthogonally general ($U$
obviously contains the open set $ X -{\bf P_{\infty}}- Sing(X)$).
\end{rem}

\begin{rem}
The situation can be slightly generalized as follows : let $Z$
be a singular projective variety, let $Z'$ be its normalization,
and let $X$ be the open set of $Z' - Sing(Z')$ where the natural
morphism to ${\bf P}^{m}$ has maximal rank. In this case, restrictions
of bundles have to be understood as pull backs. If instead one
wants to generalize to the case where $X$ is the resolution of $Z$,
many things change substantially because one does not get bundle
maps any longer.
\end{rem}

Thus we can give the following definition

\begin{df}
Let $X$ be an algebraic variety, not contained in the hyperplane
at infinity, $U$ a Zariski open set of $X$ which is orthogonally
general, and $N_{U}$ the projective normal bundle of $U$. Then
the projective normal variety $N_{X}$ of $X$ is defined
as the Zariski closure of $N_{U}$.

\end{df}

We can easily verify that the above definition is indeed independent
of the choice of $U$.

Assume now that $X$ is orthogonally general: in particular,
$X$ is smooth and we have a vector bundle (locally free sheaf)
${\cal N}_{X}$ on $X$, which is called the EUCLIDEAN NORMAL
BUNDLE of $X$.

\begin{rem}
The Euclidean Normal Bundle differs from the usual Normal Bundle
(of a smooth subvariety $ X \subset {\bf P^m}$)  defined in
algebraic geometry( cf. [Hartshorne]): the reader may in fact
notice that their respective ranks differ first of all by $1$.
However, as we shall shortly see in
the forthcoming example, they are somehow related to
each other.
\end{rem}

We can therefore compute now
the total Chern class of ${\cal N}_{X}$:\\

$c({\cal N}_{X})=c({\cal N}'_{X})$ and $c({\cal N}_{X}(-1))=
c({\cal L})c({\overline T}_{X}(-1)^{\vee})^{-1}$ by $(4)$.

 But $(3)$ yields
$c({\overline T}_{X}(-1))=c({\widehat T}_{X}(-1))$ which by $(2)$
equals $c({\cal L})c(T_{X}(-1))$.  Thus\\

 $c({\cal N}_{X}(-1))=c({\cal L})c(-{\cal L})^{-1}
c(\Omega_{X}^{1}(1))^{-1}$\\

Let us verify this formula for a hypersurface of degree $d$.  Then
we have\\

$0 \rightarrow {\cal O}_{X}(1-d) \rightarrow
\Omega_{\bf P}^{1}(1)\bigotimes {\cal O}_{X} \rightarrow
\Omega_{X}^{1}(1) \rightarrow 0$\\

and $c(\Omega_{P}^{1}(1))=c({\cal O}(1))^{-1}$.\\

So, for a hypersurface, the rank $2$ bundle ${\cal N}_{X}$ has \\

$c({\cal N}_{X}(-1))=c({\cal L})c(-{\cal L})^{-1}c(-{\cal L})
c({\cal O}_{X}(-(d-1)))=(1-H)(1-(d-1)H)$\\

( indeed, the previous formulae show ${\cal N}_{X} \cong
{\cal O}_{X} \bigoplus {\cal O}_{X}(-(d-2))$).

In general we have an exact sequence\\

$0 \rightarrow {\cal N}_{X}^{*}(1) \rightarrow
\Omega_{\bf P}^{1}(1) \bigotimes {\cal O}_{X} \rightarrow
\Omega_{X}^{1}(1) \rightarrow 0$,\\

where ${\cal N}_{X}^{*}$ is the
usual conormal bundle of $X$. \\
 Hence, $c({\cal N}_{X}(-1))=c({\cal L})c(-{\cal L})^{-1}
c({\cal N}_{X}^{*}(1))c(-{\cal L})$, and we obtain the\\

$\underline {FINAL \, FORMULA}$:
$c({\cal N}_{X}) = c({\cal N}_{X}^{*}(2))$.\\

\begin{cor}
If $X$ is a general complete intersection of degrees $d_{1}, \dots
d_{m-n}$, then $c({\cal N}_{X})= \prod_{i} (1-(d_{i}-2)H)$,
where $H$ is the hyperplane divisor.
\end{cor}

We recall once more the definition of the
Focal Locus $\Sigma_{X}$ of $X$.\\

\begin{df}
Continue to assume that $X$ is orthogonally general, let
$N_{X}\subset {\bf P}^{m}\times {\bf P}^{m}$ be the projectivization of
the Euclidean Normal Bundle, and let
$\pi =p_{2}: N_{X} \rightarrow {\bf P}^{m}$ be the second projection.
Denote then by $Y_X$ the ramification locus of $\pi$
(recall: $N_{X}$ is smooth and dim $N_{X}=m$).
 Clearly, if $X$ is projective, $Y_X \neq \emptyset$, since
$rk\, Pic (N_{X}) \geq 2$, and therefore  $\pi$ cannot be
 an isomorphism.
We define in general the focal locus as
$\Sigma_{X}:=\pi(Y_X)$.
\end{df}

\begin{df}
Let now $Z$ be any projective variety, possibly singular.
Let $X$ be a maximal orthogonally general open set of the normalization
$Z'$ of $Z$ (cf. remark 2): then the focal
locus $\Sigma_{Z}$ is defined as the closure of $\Sigma_{X}$.
\end{df}

\begin{rem}
In order to verify whether the definition would be the same
when one would replace $X$ by any orthogonally general open set
of $Z$,i.e.,
independent of the chosen open set $X$, we may observe:

\begin{itemize}
\item
For $X$ orthogonally general, the projective normal bundle $N_{X}$
has a canonical section, provided by the diagonal of $X$,
and corresponding to the tautological sheaf ${\cal L} \subset {\cal
N}_{X}(-1)$.
\item
In a neighbourhood of the canonical section, the morphism
$\pi$ is of maximal rank if and only if
 ${\cal N}'_{X}(-1) $ and $({\overline T}_{X}(-1))$ yield a
direct sum, i.e., $({\overline T}_{X}(-1))$ contains no isotropic
vectors.\\
We shall say that a point $x \in X$ is totally non isotropic if
the above situation occurs.
\end{itemize}
It follows that, in the open set of totally non isotropic points,
the ramification divisor cannot contain the fibre of the projection
 to $X$. Therefore, in this locus, the ramification divisor
 is the closure of its restriction to the inverse image of an open
set in $X$.

Instead, when there is a divisor $D$ of isotropic points of $X$,
the inverse image of $D$ may yield a component of the ramification
divisor, as happens in the following example.
\end{rem}

Consider the plane curve $C$ given, in a standard system
of Euclidean coordinates, by the parametrization $ (t, it +t^3)$.

Then the normal vector is proportional to the vector
$(i+ 3t^2, -1)$ and the endpoint map $\pi$ associates to
$(t,\lambda)$ the point

$ x = t + \lambda ( i + 3t^2)$

$y = it + t^3 - \lambda$,

and the Jacobian determinant equals

$ J = - ( 1 + 6 t \lambda) - ( i + 3t^2)^2 =
- 3 t ( 2 \lambda + 2 i t + 3 t^3)$.

Thus the focal locus consists of the evolute $\cal E$ ( image of the
curve $ \lambda = - i t - 3/2 t^3$) and of the isotropic line
$ \{ (x,y)| ix - y = 0\}$.

$\cal E$ is here the parametrical curve
$ ( 2 t - 9/2 i t^3 - 9/2 t^5, 2 i t + 5/2 t^3)$.

The previous remark and example justify the following

\begin{df}
Let now $Z$ be any projective variety, possibly singular.
We say that an  open set $X$ of $Z$ is excellent if $X$ is
 orthogonally general and $X$
is contained in the set of totally non isotropic points .
If there exists an excellent open set $X$, then the strict focal
locus $\Sigma^{s}_{Z}$ is defined as the closure of $\Sigma_{X}$.

We define instead the large focal locus $\Sigma^{L}_{Z}$
as the branch locus of the second projection $\pi$ from $N_{Z'}
\subset Z' \times {\bf P}^{m}$ to ${\bf P}^{m}$, where
as before $Z'$ is the normalization of $Z$.

Obviously one has inclusions $\Sigma^{s}_{Z} \subset \Sigma_{Z}
\subset \Sigma^{L}_{Z}$.
\end{df}

\begin{ex}
In the case of a plane curve $C$, the
strict focal locus is precisely
the evolute of the curve $C$, as in [Fantechi].
Whereas, even if all the points are totally non isotropic,
the large focal locus can be larger, as we shall now see in
the case where the curve has as a singularity a higher order cusp.
\end{ex}

Let our curve $C$ be locally given by $ (t^2, t^5)$, with respect
to some standard Euclidean coordinates; then the normal vector
is , for $t\neq 0$, proportional to  $ (-5 t^4, 2t)$, i.e., to
 $ (-5 t^3, 2)$, and thus the large focal locus is provided by
the image of the jacobian determinant of the map

$ x = t^2 -5 t^3 \lambda$,

$y =  t^5 + 2 \lambda$.

The equation of the  Jacobian determinant equals therefore

$ t ( 4 - 30 t \lambda - 25 t^6 ) =0 $, whence the large focal

locus consists of the evolute plus the line obtained for $ t=0$,
namely the $y-$ axis.

\begin{rem}
 Assume now that $Z$ is any projective variety and assume
that there is a non empty excellent open set $X \subset Z$.
If $\Sigma^{L}_{Z}$ has dimension $ \leq m-2$,
 then $\pi $ is a birational
morphism, since then $\prod_{1}({\bf P}^{m}-\Sigma^{L}_{Z})=\{1\}$.
Thus if $Z$ is not  isotropically focally degenerate and
$dim\, \Sigma^{L}_{Z}<m-1 \Rightarrow N_{X}$ is rational $\Rightarrow
Z$ is unirational, and indeed stably rational.
\end{rem}

\begin{ex} If $X$ is a smooth hypersurface of degree $d$ and
$\Sigma^{L}_{X}$ has dimension $ \leq m-2$, then $d\leq m$.\\
\end{ex}

\begin{ex}
If $X$ is a smooth complete intersection
of multidegree
$(d_{1}, \dots, d_{m-n})$, and $\Sigma^{L}_{X}$ has
dimension $ \leq m-2$, then $\sum d_{i}\leq m$.\\
\end{ex}

\begin{rem}
Let $X'\subset {\bf P}^{m}$ be a smooth variety not necessarily
satisfying the non degeneracy conditions, i.e., Assumptions $1$ and
$2$.  Then $\exists g\in {\bf PGL}(m+1)$ such that $X=gX'$ satisfies
the non degeneracy conditions.
\end{rem}
Proof\\
The non degeneracy conditions are equivalent to $(1')$ $X$ is
transversal to ${\bf P}_{\infty}$ and $(2')$ $X_{\infty}:=
X\bigcap {\bf P}_{\infty}$ is transversal to $Q_{\infty}$.  By
Bertini's theorem, we can find a hyperplane $H$ and a
smooth quadric $Q\subset H$ such that $X'$ is transversal to $H$
and $X'\bigcap H$ is transversal to $Q$.  Then choose
$h,k\in {\bf PGL}(m+1)$ such that $hH={\bf P}_{\infty}$ and
$k{\bf P}_{\infty}={\bf P}_{\infty}$, $khQ=Q_{\infty}$ and set
$g=kh$.
 $\Box$\\

Let us continue now to assume that $X$ is orthogonally general.
Moreover, we shall from now on assume that $X$ is indeed projective.
Then we can calculate
$deg\, \Sigma_{X} \, deg\, \pi _{|Y_X}$ (notice that $\pi$ is a morphism)
by
working in the Chow (or cohomology) ring of $N_{X}$.\\

Observe that, by the Leray-Hirsch theorem, the cohomology algebra
of the projective normal bundle is generated by $H^{*}(X)$
and the relative hyperplane divisor $H_2$, and holds\\

$H^{*}(N_{X})\cong H^{*}(X)[H_2]/(\sum c_{i}({\cal N}_{X}(-1))H_2^{m-n+1-i})$\\

We denote by $\Pi$ the first projection $\Pi:N_{X}\rightarrow X$, and
for commodity we also set $ p: = \pi$.

 Let
$H_{1}=\Pi^{*}(hyperplane)$, and observe that, since
${\cal N}_{X}(-1)$ is a subbundle of $ V \bigotimes {\cal O}_X$,
we have  $H_{2}=p^{*}(hyperplane)$.

 Moreover, setting $N = N_X$,  we have also the ramification formula\\

$Y=K_{N}-p^{*}(K_{\bf{P}^{m}})=K_{N}+(m+1)H_{2}$.

In order to determine the canonical divisor $K_{N}$ of $N =  N_X$,
we write as
usual\\

$K_{N}=K_{N|X}+ \Pi^{*}(K_{X})$,\\

where $K_{N|X}$ can be calculated through the Euler exact sequence
for the relative tangent bundle $T_{N|X}$ of $N$\\

$0 \rightarrow {\cal O}_{N}(-H_2) \rightarrow \Pi^{*}({\cal N}_{X}(-1))
\rightarrow T_{N|X}(-H_2) \rightarrow 0$,\\

whence\\

$K_{N|X}=-c_{1}(T_{N|X})=-[c_{1}(T_{N|X}(-H_2)+(m-n)H_2)]=
-(m-n+1)H_2 -\Pi^{*}(c_{1}({\cal N}_{X}(-1)))= -(m-n+1)H_2+H_{1}-
\Pi^{*}(c_{1}({\cal N}_{X}^{*}(1)))=
-(m-n+1)H_2 +H_{1}-c_{1}({\cal N}_{X}^{*})-(m-n)H_{1}=
-(m-n+1)H_2 +H_{1}+(m+1)H_{1}+K_{X}-(m-n)H_{1}$\\

In the end we obtain:\\

$K_{N}=2\Pi^{*}K_{X}-(m-n+1)H_2 +(n+2)H_{1}$\\ thus
we get the \\

$\underline{CLASS-FORMULA}$: $Y_X =2\Pi^{*}K_{X}+nH_2 +(n+2)H_{1}$,\\

and the\\

$\underline{DEGREE-FORMULA}$:\\
\\
  $( deg\, \Sigma \,) ( deg\, p_{|Y}) =
H_2^{m-1}(2\Pi^{*}K_{X}+nH_2 +(n+2)H_{1})$.\\

In the sequel ( section 5) we shall see how the above cited
Leray-Hirsch Theorem
allows to make the degree formula more explicit.

\section{Non degeneracy of Focal Loci}

Throughout this section we assume that $X$ is projective
and orthogonally general, i.e.,  the non
degeneracy conditions  $0-2$ above are satisfied, in
particular we have that
 $N=N_{X}$ is a  bundle .  Our aim is then to
determine for which $X$ it is possible that $\Sigma_{X}$ is
degenerate, that is, has dimension strictly less than $m-1$.
It is easy to see that, if $X$ is a linear space, then $\Sigma_{X}$
 is degenerate
and is a linear space of dimension equal to $codim\, X-1$.  In
what follows, we shall prove that if $X$ is orthogonally general
also the converse holds, i.e., if
$\Sigma_{X}$ is strongly degenerate, then $X$ is a linear space.\\

\begin{nota}
Let ${\bf C^{m}}={\bf P^{m}}-{\bf P}_{\infty}$, $N_{\infty}=
p^{-1}({\bf P}_{\infty})$, $N_{a}=p^{-1}({\bf C^{m}})$.
\end{nota}

We have

\begin{lem}
$dim\,p^{-1}(y)=0$ $\forall y \in {\bf C^{m}}$
\end{lem}
Proof\\
After identifying $p^{-1}(y)$ with the set
$\Gamma=\{x\in X: y\in N_{x}\}$, it is easy to see that $\Gamma$
has empty intersection with the hyperplane ${\bf P}_{\infty}$.
Indeed, if  $x\in X_{\infty}$, then
$N_{x}\subset {\bf P}_{\infty}$   $\Box$\\

\begin{cor}
If $\Sigma$ is a component of the focal locus,
image of a component $Y$ of the divisor $Y_X$,
and $dim\, \Sigma<m-1$, then\\
$(i)$ $Y\subset N_{\infty}$
(since $\forall y\in\Sigma$ dim $p^{-1}(y)>0$)\\
$(ii)$ $\Sigma \subset \neq {\bf P}_{\infty}$ (hence $\Sigma$ is
degenerate).\\
$(iii)$ if  for every component $\Sigma$   of the focal locus
holds  $dim\, \Sigma<m-1$, then\\  $p:N_{a}\rightarrow {\bf C^{m}}$
is an isomorphism.
\end{cor}

\begin{rem}
The divisor $N_{\infty}$ splits as $N_{|X_{\infty}}\bigcup N'$,
where $N'={\bf P}({\cal N}_X')$.
\end{rem}

Let us first consider the case where $X$ is a curve
( for this case we shall give a different proof in the sequel,
showing that then either $C$ is a line, or $C$ is a circle, what
contradicts the hypothesis that $C$ be orthogonally general).

CASE: $X=C$ curve.\\

Let $C$ be an irreducible ( and orthogonally general)
 curve of degree $d$.
Then $N_{|C_{\infty}}$ consists of
$d$ distinct copies of ${\bf P}_{\infty}$,
$p:N_{|C_{\infty}}\rightarrow {\bf P}_{\infty}$
is a finite map, and by the transversality of $C$ to ${\bf P}_{\infty}$,
the divisor $Y_C$ does not contain any component of $N_{|C_{\infty}}$.

  Therefore we get\\

\begin{cor}
If $C$ is an irreducible ( and orthogonally general)
curve and $dim\, \Sigma_C <m-1$, then $Y_C =N'$
(set-theoretically)
\end{cor}
Proof\\
Indeed, $Y_C \subset N_{\infty}$, but no component of $N_{|C_{\infty}}$
 is contained in $Y_C$. Thus $Y_C \subset N'$.  We can conclude since
 $dim\, Y_C =dim\, N'$ and $N'$ is irreducible (being
a projective bundle on the curve $C$).
$\Box$\\

\begin{prop}
Assume again that $C$ is an irreducible ( and orthogonally general)
curve . Then
 $(1)$ $dim\, \Sigma<m-1 \Leftrightarrow $ $(2)$ $ C$ is a line.\\
 $(3)$ In this case $\Sigma$ is a linear space of dimension
$m-2=codim\, C-1.$
\end{prop}
Proof\\
(1) $\Leftarrow$  (2) being clear, let's prove the other implication
$(1)$ $\Rightarrow$ (2): \\
 Let $N'_{p}$ be the fibre of $N'$ over
$p\in C$, which is a hyperplane in ${\bf P}_{\infty}$.  Now
$\Sigma_C =p(N')$ is irreducible, has dimension $<m-1$ and
contains $N'_{p}$, which has dimension equal to $m-2$.  Therefore,
 $\Sigma_C =N'_{p}$ and  $N'_{p}=N'_{q}$
$\forall p,q\in C$.  This implies
$T_{p}C\bigcap {\bf P}_{\infty}=T_{q}C\bigcap {\bf P}_{\infty}$
$\forall p,q \in C$.

This clearly implies that $C$ is a line, since then
for each point $p \in C$ the projective tangent line $T_{p}C$
is the join of $p$ and of a fixed point $p_{\infty}$
(whence one can find then
$m-1$ independent linear forms vanishing on $C$).

  $\Box$\\

CASE: $dim\, X=n\geq 2$\\
Since $X_{\infty}$ is smooth, by Bertini's theorem $X_{\infty}$ is
irreducible. Therefore also $N|_{X_{\infty}}$
and $N'$ are irreducible.\\
We have

\begin{lem}
If $ n \geq 2$ and  $dim\, \Sigma_X < m-1$, then\\
(i) $Y_X =N'$ set-theoretically.\\
(ii) $[Y_X]=[nN']$ in $Pic(N)$.\\
(iii) $2(K_{X}+(n+1)H)=0$ in $Pic(X)$, where $H = H_1$ is
the hyperplane
divisor on $X$.
\end{lem}

Proof\\
Since $p$ is surjective, we have one and only
one of the following two
cases:
(a) $p(N|_{X_{\infty}})\subset\neq {\bf P}_{\infty}$; (b)
$p(N')\subset\neq {\bf P}_{\infty}$. But (a) cannot hold.
Indeed,
since $[H_{1}]=[N|_{X_{\infty}}]$ in $Pic(N)$, (a) implies
$[Y]=[\alpha H_{1}]$  for some $\alpha >0$.  But then from the class
 formula $(*)$ $Y_X =2\Pi^{*}K_{X}+nH_2 +(n+2)H_{1}$, it follows that
$2\Pi^{*}K_{X}+ nH_2 +(n+2- \alpha) H_{1} = 0$, contradicting
the Leray-Hirsch Theorem.

Therefore,
(b) holds and hence $Y=N'$ set-theoretically. (ii) and (iii)
follow immediately from the class formula
$(*)$, because $H_2 =H_{1}+[N']$ in $Pic(N)$.  $\Box$\\

>From point $(iii)$ it follows that

\begin{cor}
If $dim\, \Sigma_X< m-1$, then $X$ is a linear space.
\end{cor}
Proof\\
Let $C=X\bigcap H_{1}\bigcap \dots \bigcap H_{n-1}$ be a smooth curve.
By successive applications of the adjunction formula $(iii)$ yields
$2(K_{C}+2H)=0$.  Extracting degrees, we get $2(2g(C)-2+2\, deg(C))=0$,
which is equivalent to $g(C)=0$ and $deg(C)=1$.  $\Box$\\

We can conclude

\begin{teo}\label{orthgen}
Let $X\subset {\bf P}^{m}$ be a projective variety of dimension
$n \geq 1$ which is orthogonally general.  Then
$dim\, \Sigma_{X}<m-1 \Leftrightarrow X$ is a linear space.
In this case, $\Sigma_X$ is a linear space of dimension equal to
$codim\, X-1$.
\end{teo}

\section{The Degree of the Focal Locus of a Surface}
Let $X^{2}=S\subset {\bf P}^{m}$ be a surface and assume that
$S$ satisfies the non-degeneracy conditions. Setting $n=2$ in the
Degree-Formula given in Section $1$, we get ( recall $ H = H_1$)\\

 $(F1)$ $deg\, \Sigma_S \, deg\, p_{|Y_S}=2H_2^{m-1}(K_{S}+H_2+2H)$\\

Our first aim in this section is to express the right-hand side of
$(F1)$ in terms of the Chern classes $c_{1}(S)$, $c_{2}(S)$ and
of the hyperplane divisor $H$ of $S$.\\

By the Leray-Hirsch theorem $H_2^{m-1}=-c_{1}({\cal N}_{S}(-1))H_2^{m-2}
-c_{2}({\cal N}_{S}(-1))H_2^{m-3}$.

 Using this relation, the right-hand
side of $(F1)$ becomes\\

$(*)$ $2H_2^{m-2}(c_{1}({\cal N}_{S}(-1))^{2}-c_{2}({\cal N}_{S}(-1))-
K_{S}c_{1}({\cal N}_{S}(-1))-2H c_{1}({\cal N}_{S}(-1)))$\\

Recall that $c({\cal N}_{S}(-1))=c({\cal L})c({\cal N}^{*}_{S}(1))$, where
${\cal L} ={\cal O}_{S}(-1)$ and ${\cal N}^{*}_{S}$ is the conormal bundle
of $S$. Thus,\\

$(1)\begin{array}{l}
c_{1}({\cal N}_{S}(-1))=c_{1}({\cal N}^{*}_{S}(1))-H=
c_{1}({\cal N}^{*}_{S})+(m-3)H\\
c_{2}({\cal N}_{S}(-1))=c_{2}({\cal N}^{*}_{S}(1))-
H c_{1}({\cal N}^{*}_{S}(1))=\\
c_{2}({\cal N}^{*}_{S})+(m-4)Hc_{1}({\cal N}^{*}_{S})+
{\frac{1}{2}}(m-2)(m-5)H^{2}
\end{array}$\\

Using the normal-bundle sequence we get\\

$(2)\begin{array}{l}
c_{1}({\cal N}^{*}_{S})=c_{1}(S)-(m+1)H\\
c_{2}({\cal N}^{*}_{S})=-c_{2}(S)+{\frac{1}{2}}m(m+1)H^{2}+
c_{1}(S)c_{1}({\cal N}^{*}_{S})=\\
-c_{2}(S)+{\frac{1}{2}}m(m+1)H^{2}+c_{1}(S)^{2}-(m+1)Hc_{1}(S)
\end{array}$\\

and substituting in $(1)$, we get\\

$(3)\begin{array}{l}
c_{1}({\cal N}_{S}(-1))=-4H+c_{1}(S)\\
c_{2}({\cal N}_{S}(-1))=9H^{2}-5Hc_{1}(S)+c_{1}^{2}(S)-c_{2}(S)
\end{array}$\\

Hence $(*)$ becomes\\

$(**)$ $2H_2^{m-2}(15H^{2}+c_{1}^{2}(S)+c_{2}(S)-9Hc_{1}(S))$\\

We recall that $H_2 =[N']+H$, so that $(**)$ can be rewritten as\\

$(***)$ $2[N']^{m-2}(15H^{2}+c_{1}^{2}(S)+c_{2}(S)-9Hc_{1}(S))$\\

Finally, since $[N']^{m-2}_{|F}=1$, where $F$ is a fibre of
>$\pi:N\rightarrow S$, we conclude\\

$(DF)$  $deg\, \Sigma_S \, deg\, p_{|Y_S}=
2(15H^{2}+c_{1}^{2}(S)+c_{2}(S)-9Hc_{1}(S))=
2(15d+c_{1}^{2}(S)+c_{2}(S)-9Hc_{1}(S))$,\\

where $d=deg(S)$.\\

By Noether's formula, we can also write\\

$(DF')$  $deg\, \Sigma_S \, deg\, p_{|Y_S}=
2(15d+12{\chi}({\cal O}_{S})-9Hc_{1}(S))$.\\

We can express also our formula in terms of the  sectional genus
$ \pi$ of our surface $S$ ( recall that $2 \pi - 2 =
H^2 - Hc_{1}(S)$) as

$(DF'')$  $deg\, \Sigma_S \, deg\, p_{|Y_S}=
2( 18(\pi - 1) + 6d+12{\chi}({\cal O}_{S}))$\\

\begin{ex}
  For $m=3$, we have
 $c_{1}(S)=(4-d)H$ and

$c_{2}(S)=6H^{2}-d(4-d)H^{2}$.  Thus\\

$deg\, \Sigma_S \, deg\, p_{|Y_S}=2d(d-1)(2d-1)$\\
\end{ex}

\begin{ex}
 For $m=4$, we have the formula
$c_{2}(S)=c_{1}(S)^{2}-5Hc_{1}(S)+10d-d^{2}$ $[Hartshorne, p. 434]$,
or , equivalently,

$d^2 - 5 d + 2 ( 6 {\chi}({\cal O}_{S}) -c_{1}(S)^{2}) =
10 ( \pi -1)$

which gives\\

$deg\, \Sigma_S \, deg\, p_{|Y_S}=2/5 ( 9 d^2 - 15 d +
168 {\chi}({\cal O}_{S}) - 18  c_{1}(S)^{2})$\\

\end{ex}

\section{Weakly focally degenerate varieties}

In this section we shall first consider the case of a hypersurface
$X$ of dimension $n$ , and we shall characterize the case
where $X$ is weakly focally degenerate.
The characterization of the hypersurfaces $X$ which are
 isotropically focally degenerate
  will be given in the next section.

Later on in this section we shall deal with focally degenerate
varieties of any codimension.

 We shall essentially use very classical tools such as
  the implicit function theorem, dimension counts and the
standard method of obtaining new equations by differentiating
old ones .

Let $F(x_1,... x_{n+1}) = 0$ be the affine polynomial equation
of a hypersurface $X$. We shall in this section be mostly
interested about a birational description of $X$, whenceforth
we might, by abuse of notation, not distinguish between a projective
variety and its affine part (or any nonempty Zariski open set of it).

In this case the gradient $\nabla F$ of $F$ gives a trivialization
of the Normal Bundle $N_X$ at the smooth points of $X$,
 and the second projection $\pi:N_{X}
 \rightarrow {\bf P^{n+1}}$
coincides with the endpoint map \\

$\epsilon (x,\lambda) = x + \lambda \nabla F(x)$,

where $x= (x_1,... x_{n+1}) $ is a point of $X$, and
$\lambda$ is a scalar coordinate $ = \lambda_{1} /\lambda_{0}  $ ,
$ (\lambda_{0} ,\lambda_{1})  $ being homogeneous coordinates
on ${\bf P^{1}}$.

As a warm up, let us investigate when does it occur that
the endpoint map is not finite.  That is, let us assume that
$\Gamma$ is a curve in $N_X$ which is mapped to
a point $O$ by the endpoint map $\epsilon$, and that this point
does not lie at infinity.

Choosing a parameter $t$ for $\Gamma$, we have functions
$ x(t), \lambda(t)$ such that

1) $ F (x(t)) \equiv 0 $\\
1')$ x(t) + \lambda(t) \nabla F ( x(t)) \equiv O  $.

If $x(t)$ is a smooth point of $X$, then the gradient
$ \nabla F ( x(t)) $ does not vanish, whence $x(t)$ is not constant:
thus at a general point of $\Gamma$ we may assume that the
derivative $ \dot x(t) : = d x(t)/dt$ is non vanishing.

Let us use the scalar product $ < , >$ associated to the
quadratic form $Q_{\infty}$, and let us choose affine
coordinates such that $ < , >$ is the standard scalar
product ( i.e., the matrix of $Q_{\infty}$ is the
identity matrix); since

2) $ x(t) - O \equiv - \lambda(t) \nabla F ( x(t))  $ , and
$ < \nabla F ( x(t)) , \dot x(t) > \equiv 0 $ we infer that

3) $< x(t) - O , x(t) - O > \equiv constant$.

Therefore, the basis curve $ \gamma \subset X = \{ x| F(x) =0\}$
is a curve contained in a sphere with centre the point $O$
( note that the sphere may also have radius zero !).

Conversely, if we have such a spherical curve $\gamma$ meeting $X$
and with the property
that the two vectors  $ x(t) - O , \nabla F ( x(t))  $ are
proportional, then we find $\lambda(t)$ so that 1'), 1) hold,
whence we find $\Gamma$ which is mapped to the point $O$ by
the endpoint map ( and moreover it follows from 1) that
$\gamma$ is contained in $X$).
Finally, since $\Gamma$ is mapped to a point, it is obviously
contained in the ramification divisor $Y$ of the endpoint map.

We have therefore the following

\begin{prop}
Given a smooth affine hypersurface $X$ , the
positive dimensional irreducible components of
the  fibres $\pi ^{-1} (O)$
of the map to the affine part of the Focal Locus correspond exactly to the
subvarieties $\Phi$
contained in a sphere $S$ with centre $O$, and such that
$X$ is everywhere tangent to $S$ along $\Phi$.

\end{prop}

Proof\\
Let $\Psi$ be a component of the fibre $\pi ^{-1} (O)$.
Then consider that $\Psi$ is the union of the curves $\Gamma$
contained in it : each of these projects to $\gamma \subset X$
contained in a sphere $S_c$ with centre $O$ and radius $c$.
But the image of $\Psi$, call it $\Phi$, is irreducible, whence
all the radii are equal, and we get the desired sphere $S$.
Conversely, the tangency condition provides a rational function
$\lambda$ on $\Phi$ whose graph is the required variety $\Psi$.
 $\Box$\\

It is rather clear that the previous proposition allows easily
to construct examples where the map $\pi : Y \rightarrow \Sigma_X$
is not finite.

\begin{rem}
If instead the point $O$ is at infinity, let's identify it with
one vector in $V'$, then we get the equation

 $ O \equiv  \lambda(t) \nabla F ( x(t))  $, whence

$ < O ,  x(t) > \equiv constant $ .

So, in this case, the
positive dimensional irreducible components of
the  fibres $\pi ^{-1} (O)$ correspond exactly to the
subvarieties $\Phi$
contained in a hyperplane $H$ with normal direction $O$, and such that
$X$ is everywhere tangent to $H$ along $\Phi$.

\end{rem}

We can push the previous calculations to describe the
weakly focally degenerate hypersurfaces.

Let us thus assume that $X =  \{ x| F(x) =0\}$ is weakly
focally degenerate. This simply means that there is a component
$\Sigma$ of the focal locus which has dimension

$  dim \Sigma = a < n$.

Arguing as before, we notice that $\Sigma$ will simply be any
maximal irreducible variety such that its inverse image
in $N_X$ has a
dominating component $Z$ of dimension  $n$.
We can analogously treat the case where this dimension is bigger than
$n$, i.e., when $Z$ = $N_X$, or equivalently  $X$ is
isotropically focally degenerate : in this case we may also allow
$  dim \Sigma = n$.

We have thus an irreducible component $Z$ of the ramification
divisor, with $\pi (Z) = \Sigma$.

To start with, let us assume that $ \Sigma \not\subset {\bf P_\infty}$.

Therefore, at the general point of $Z$ we can choose local
coordinates

$s=(s_1, ...., s_a)$ and  $ t= (t_1, .. t_{\nu-a})$  ( $\nu = n $ or $n+1$)\\

such that the fibres of $\pi$ are locally given by setting
$ s = constant$, in other words we have functions

$ x(s,t) , \lambda(s,t) $ parametrizing the points of $Z$ ,

and a function $ O(s)$ parametrizing the image $\pi (Z) = \Sigma$
of the end-point map. This means that the following
equations hold :

1") $ F (x(s,t)) \equiv 0 $\\

2") $ x(s,t) - O(s) \equiv - \lambda(s,t) \nabla F ( x(s,t))  $ ,

differentiating 1") with respect to both sets of variables $s,t$,
we infer that

$ < \nabla F ( x(s,t)) ,( d x(s,t)/ d t_i )> \equiv 0 $ as well as

$ < \nabla F ( x(s,t)) ,( d x(s,t)/ d s_j) > \equiv 0 $.

We argue as we did before :

since $ x(s,t) - O(s)$ is proportional to $\nabla F ( x(s,t))  $,
we obtain that $ x(s,t) - O(s)$ is orthogonal to all the partial
derivatives of $x(s,t)$.

Since however $( d x(s,t)/ d t_i ) = ( d (x(s,t) - O(s))/ d t_i )$,
it follows that there is a function $r(s)$ such that

3") $< x(s,t) - O(s) , x(s,t) - O(s) > \equiv r(s)$.

What we have done insofar is to write down the family of spheres
containing the projections $X_s$ to $X$ of the fibres over
$O(s) \in  \Sigma$.

On the other hand, we can use the other partial derivatives
$( d x(s,t)/ d s_j)$ in order to obtain a complete description
of $X_s$.

In fact, let us calculate the partial derivatives $(\partial r(s)/
 \partial s_j)$

 They are $= 2 < x(s,t) - O(s), ( \partial (x(s,t) - O(s)/
\partial s_j)>$

$= - 2 < x(s,t) - O(s), ( \partial O(s)/ \partial s_j)>$.

We have therefore established

4") $(\partial r(s)/ \partial s_j) = - 2 < x(s,t) - O(s), ( \partial
 O(s)/ \partial s_j)>$,

whose geometric meaning is the following: if $O(s)$
is a smooth point of $\Sigma$, whence all the partial derivatives
$( d O(s)/ d s_j)$ are linearly independent, then
 $X_s$ is contained in the intersection of the sphere given
by 3") with the codimension $a$ affine subspace given
by 4").

If this intersection has the expected dimension $ n-a$, then
it has the same dimension as $X_s$ and if it is moreover
irreducible it will coincide with $X_s$.

\begin{lem}
Consider an affine subspace $ L = \{ x | < x- O, v_j> = c_j$
for $ j=1,.. a \}$ of codimension $a$ and assume that
$L$ is contained in the sphere $S(O,r^{1/2}) =
  \{ x | < x- O, x-O>= r \}$.
Then

(*) the direction $W$ of $L$ is an isotropic subspace for $< , >$,
and there exists $ x_0 \in L$  such that $ x_0 - O$
is orthogonal to $W$ ( equivalently, $W$ is isotropic and
$ L \subset O + W^{\perp}$).

Observe
moreover that the orthogonal $W^{\perp}$ is the vector  space
$U$ generated by the vectors $v_j$.

Also the converse holds, in the sense that
   if (*) is verified, then there exists a constant $R$
such that $L$ is contained
in the sphere $S(O,R^{1/2})$.

\end{lem}
Proof\\
Let $ x_0 \in L$  and write $ L = x_0 + W$.
Since $ <x-O, x-O> \equiv r$ for $ x \in L$, we get

$ <x_0-O, x_0 -O> + 2 < w,x_0 -O> + < w, w> \equiv r$
for each vector $ w \in W$.

Thus the quadratic polynomial $ <w,w>$ is identically
zero on $W$ , what amounts to say that the subspace
$W$ is isotropic; the vanishing of the linear form
yields the desired orthogonality of $x_0 - O$ to $W$.

Conversely , $ <x-O, x-O> \equiv  <x_0-O, x_0 -O> $
and $L$ is contained in the sphere

$  \{ x | < x- O, x-O>= R \}$ once we set
 $ R =<x_0-O, x_0 -O> $.
$\Box$\\

\begin{lem}
Consider an affine subspace $ L = \{ x | < x- O, v_j> = c_j\}$
as in the previous lemma $3$, and assume that the affine quadric
$L \cap S(O,r^{1/2})$ is reducible.
Then either

(i)  $dim (W/ W \cap W^{\perp} ) = 1$
and there exists $ x_0 \in L$  such that $ x_0 - O$
is orthogonal to $W$ and  $ < x_0- O, x_0- O> \neq r$
 or

(ii)  $dim (W/ W \cap W^{\perp} ) = 2$
and there exists $ x_0 \in L$  such that $ x_0 - O
\in W^{\perp}$ and  $ < x_0- O, x_0- O> = r$

\end{lem}
Proof\\
As before , for each choice of $ x_0 \in L$  we can write $ L = x_0 + W$.
Since the equation of our affine quadric is

$ <x_0-O, x_0 -O> + 2 < w,x_0 -O> + < w, w> = r$
for each vector $ w \in W$,
and we impose the condition that the quadric be the union of two
affine hyperplanes, it follows
that the quadratic form $ <w,w>$ on $W$ has  rank
either $1$ or $2$.

In the latter case,
since the rank of the complete quadric equals the rank
of the quadratic form, acting with a translation on $W$,
  we can kill the terms of lower degree.

In the former case, if the linear part of the equation would
not belong to the image under $Q_{\infty}$ of $W/ W \cap W^{\perp}$,
the rank would be at least $3$. Whence, acting with a
 translation on $W$, we may kill the linear part and then the constant
must be non zero.

$\Box$\\

We have therefore found that the projection of $Z$
is contained in the locus $X'$ given by

3"") $ \{x|\exists s, < x - O(s) , x - O(s) > \equiv r(s)$

4"") $(d r(s)/ d s_j) = - 2 < x - O(s), ( \partial
 O(s)/ \partial s_j)>\}$.

If moreover $Z$ surjects onto $X$ and $X'$ is irreducible, then
 $X'$ equals $X$ unless we are in the
exceptional case where (cf. Lemma 3) for each point
$O(s)$ the (vector) tangent space
$V_s$ to $\Sigma$ at $O(s)$ satisfies the condition that
$V_s$ contains its orthogonal $W_s : = V_s^{\perp}$, and moreover
then $( x(s,t) - O(s) )$, for each $t$ belongs to the
subspace $V_s: = W_s^{\perp}$.

The locus $X'$, as written, is the projection of the locus

$Z' \subset {\bf P}^{m} \times \Sigma $ defined as

3"') $ \{(x,s)| < x - O(s) , x - O(s) > \equiv r(s)$

4"') $(\partial r(s)/ d s_j) = - 2 < x - O(s), ( \partial
 O(s)/ \partial s_j)>\}$.

If we calculate the tangent space to $Z'$ at the point $(x,s)$
we obtain that it is contained in the hyperplane:

5"') $\{(\xi,\sigma)| 2 < x - O(s) ,\xi > -2
< x - O(s) ,( \partial O(s)/ \partial \sigma) > - ( \partial
 r(s)/ \partial \sigma) = 0\}$
=  $\{(\xi,\sigma)| 2 < x - O(s) ,\xi >  = 0\}$

(since $(x,s)$ is a point of $Z'$).

By Sard's lemma, $X'$ has dimension at most $n$: whence,
if we assume that the component $Z$ dominates $X$,
and thus $X \subset X'$, we conclude that $X= X'$
( in the  exceptional case, or if $X'$
is irreducible) or at least that $X$ is a component of $X'$.

We are now in the position to explain the main
constructions which are underlying the characterization of
 the focally  degenerate  varieties.

\begin{df}

THE INVERSE CONSTRUCTION TO  FOCAL DEGENERACY.

Start from  the following data :

i) Let $\Sigma$ be an irreducible affine variety of dimension
$a$ , and let $\Sigma^{\star}$ be an irreducible subvariety
of the product $\Sigma \times \bf C$ which is the graph
of an algebraic function $r$ on $\Sigma$.

Proceed constructing an algebraic set $X'$ as follows:

ii) The subvariety $\Sigma^{\star}$ defines a family of spheres

$Z^{\star} \subset {\bf C}^{m} \times \Sigma \times {\bf C}$
defined as

3"') $ \{(x,O,r)|(O,r) \in \Sigma^{\star},
 < x - O , x - O > = r\}$.

iii) Consider in $ {\bf C}^{m} \times  {\bf C}$ the tangent
space to $\Sigma^{\star}$ at a point $(O,r)$, and its
orthogonal with respect to the quadratic form $Q_\infty
\bigoplus 1$: under the embedding of $ {\bf C}^{m}$ in
$ {\bf C}^{m} \times  {\bf C}$ sending $x$ to
$(x-O,1/2)$, its pull back is precisely an affine space
given by an equation as 4"').
We can in this way define a bundle (if $\Sigma$ is smooth)
of affine spaces

$A^{\star} \subset {\bf C}^{m} \times \Sigma^{\star}$,

$A^{\star} = \{(x,O,r)|(O,r) \in \Sigma^{\star},
 ( x - O , 1/2)  \in {T \Sigma^{\star}_{(O,r)}}^{\perp} \}$.

iv) define $Z'$ as the intersection $Z^{\star}  \cap A^{\star}$
( a divisor in $A^{\star}$)

v) define $X'$ as the projection of $Z'$ on the first factor
${\bf C}^{m}$;

vi) observe that, by the argument we gave above,
$ dim X' \leq m-1$.

vii) assume finally that  $\Sigma , r$ are $\bf admissible$,
which amounts to the requirement that $Z'$ dominate  $\Sigma$.
\end{df}

\begin{rem}
The condition that $\Sigma , r$ be $\bf admissible$ is obviously
satisfied unless $Z'$ is a union of fibres of the projection
$A^{\star} \rightarrow \Sigma$. This means, unless the quadratic
function $ < x - O(s) , x - O(s) > $ is constant on the affine spaces
$A^{\star} _s$. Therefore, the pair  $\Sigma , r$ is $\bf admissible$
unless we are in the situation of Lemma 3 , whence
$ < x - O(s) , x - O(s) > \equiv R(s) $  on $A^{\star} _s$,  but
$ R(s) \not\equiv r(s)$.
\end{rem}

There remains however to see what happens in the case
 where $\Sigma$ lies at
infinity .

In this case, we derive (cf. remark 8) the following equations,
where $O(s)$ is a $V'$ -valued function leading to a parametrization
of $\Sigma$ :

6)  $  < x (s,t) , O(s) > \equiv r(s)$

7) $(\partial r(s)/ \partial s_j) \equiv < x(s,t) , ( \partial
 O(s)/ \partial s_j)>, for  j = 1,... a$.\\

In this case, if $O(s)$ is a smooth point of $\Sigma$, then
the  $a+1$ vectors $ O(s) ,\partial
 O(s)/ \partial s_j$ are linearly independent and 6) and 7)
imply that $X_s$ is contained in the affine space

8) $ X'_s = \{ x| < x , O(s) > \equiv r(s) ,\\
(\partial r(s)/ \partial s_j) \equiv < x , ( \partial
 O(s)/ \partial s_j)>, for  j = 1,... a \}$.

Since $\Sigma$ lies at infinity , $X$ is not isotropically
focally degenerate, whence $Z$ has dimension $m-1$ : it follows
that $X_s$, $X'_s$ have the same dimension $m-1-a$, whence
they coincide.

Moreover, $Z$ must dominate $X$, else a whole fibre of
$N_X \rightarrow X$ is contained in $Z$, and therefore
its projection cannot lie at infinity (remember that $X$ is here
supposed to be affine).

Therefore, it follows that $X$ equals $X'$, the closure of the
union of the $X'_s$.

We are therefore led to the following
\begin{df}

THE ASYMPTOTIC INVERSE CONSTRUCTION TO WEAK FOCAL DEGENERACY.

Start from  the following data :

i) Let $\Sigma$ be an irreducible  variety of dimension
$a$ , contained in ${\bf P_{\infty}}$, and let $\Sigma^{\star}$
 be an irreducible subvariety
of the product $\Sigma \times \bf C$ which is the graph
of an algebraic section $r$ of ${\cal O}_{\Sigma}(1)$.

Then consider the algebraic set $X'$ which is the closure
of the union of the family of affine spaces $X'_s$
defined by  8).

\end{df}

\begin{rem}
The attentive reader will find a slight abuse of notation above,
which can be explained as follows : in the case where $\Sigma$
does not lie at infinity, since we have a privileged affine
chart = ${\bf P^m}- {\bf P_{\infty}}$, we consider $r(s)$ just as
an algebraic function on $\Sigma$. If however $\Sigma \subset
{\bf P_{\infty}}$, then there is no favourite standard affine chart
and we  make clear that $r$ is not really a function, but
a  section of ${\cal O}_{\Sigma}(1)$
 (possibly multivalued and with poles!).
\end{rem}

\begin{rem}
Consider a variety $X=X'$ obtained from the asymptotic inverse
focal construction.

Then its part $ X_{\infty} = X \cap
{\bf P_{\infty}}$ consists of the points \\

$ \{ x \in {\bf P_{\infty}}| < x , O(s) > \equiv r(s) ,\\
(\partial r(s)/ \partial s_j) = < x , ( \partial
 O(s)/ \partial s_j)>, for  j = 1,... a \}$.

If we therefore identify ${\bf P_{\infty}}$ with its dual space
via the quadratic form $Q_{\infty}$, it follows that $X_{\infty}$
is the dual variety of $\Sigma$ !

Observe moreover, that if $X$ is a linear space, then $\Sigma_X$
 equals  $X_{\infty}^{\perp}$.

In this case the section $r(s)$ is just induced by  a linear form
on $\Sigma_X$ (i.e., a vector in $ (V')^{\vee}$).
\end{rem}

We observe now that we have  insofar proved the following

\begin{teo}
Let $X$ be a  focally degenerate hypersurface in
${\bf C}^{n+1}$ and let $\Sigma$ be a component of the strict focal
locus  ( i.e., we are in the non vertical case and the corresponding
component $Z$ of $Y_X$ projects onto $X$).
 Then
\begin{itemize}
\item
either $\Sigma$ is contained in ${\bf P_{\infty}}$ and
 $X$ is obtained from $\Sigma, r$ via the asymptotic
inverse focal construction associated to an
algebraic section $r$ of ${\cal O}_{\Sigma}(1)$
\item
or, $\Sigma$ is not contained in ${\bf P_{\infty}}$ and
there is an algebraic function $r(s)$ on $\Sigma$ such that,
 applying the inverse construction to focal degeneracy,
we get a hypersurface $X'$ such that  $X$ is a component of $X'$.
\end{itemize}

Conversely, start from any  admissible pair $( \Sigma, r)$,
 and assume that
an irreducible hypersurface is a component of the algebraic
set $X'$ obtained from the inverse construction or from
the asymptotic inverse construction : then
$X$ is a  focally degenerate hypersurface.
\end{teo}

Proof\\
There remains only to show that if $X$ is
an irreducible hypersurface, component of the algebraic
set $X'$ obtained from an inverse construction : then
$\Sigma$ is a component of the focal locus of $X$.
This follows since, by 5"'), $ x- O(s)$ is a normal
vector to $X'$, respectively since $O(s)$ is a normal vector
to $X'$; moreover, $Z'$ dominates $\Sigma$ by the assumption
that $r$ be admissible.
$\Box$\\

However, the inverse constructions, as we are going to see,
work more generally also in the case where $X'$ has smaller
dimension than the expected dimension $m-1$.

We have in fact the following

\begin{teo}
Let $X$ be a  focally degenerate variety of dimension $n$ in
${\bf C}^{m}$ and let $\Sigma$ be a component of the focal
locus of dimension $a \leq m-1$, projection of a
component $Z$ of $Y_X$.
Then $\Sigma$ determines birationally an
irreducible subvariety
$ \Sigma^{\star}$ of $\Sigma \times \bf C$ corresponding to an
algebraic section $r(s)$ of ${\cal O}_{\Sigma}(1)$ and,
applying the appropriate inverse construction to focal degeneracy,
we get  an algebraic set $X'$ which is focally degenerate,
and  indeed isotropically focally degenerate in the case where
 $\Sigma$ is
not contained in the hyperplane at infinity ${\bf P_{\infty}}$
and $dim Z'= m$ (in this case the fibres
$X'_s$ of $Z' \rightarrow \Sigma$ are affine spaces).

 There are seven cases :
\begin{itemize}
\item
1) $X$ is  isotropically focally degenerate : then $X=X'$, $dim Z'= m$ and
the fibres $X_s$ of $N_X \rightarrow \Sigma$ are affine spaces.
Moreover, here $\Sigma$ is not contained in ${\bf P_{\infty}}$.
\item
2) $\Sigma$ is not contained in ${\bf P_{\infty}}$ ,
$Z$ projects onto $X$ and $X'$ is not isotropically focally degenerate:
then $X$ is a component of $X'$
\item
3) $Z$ projects onto $X$, $X'$ is  isotropically focally degenerate,
but $X$ is not isotropically focally degenerate:
then  $X \subset X'$ is a divisor, $Z$ is the restriction to
 $X$ of the normal bundle $N_{X'}$,
and $\Sigma$ is the focal locus of $X'$ (again here
$\Sigma$ is not contained in ${\bf P_{\infty}}$)
\item
4) $\Sigma$ is not contained in ${\bf P_{\infty}}$ ,
$Z$ projects onto a divisor $X"\subset X$, $X"$ is a component
of $X'$, $X"$ is focally
degenerate, with a component $Z"$ of the ramification locus $Y_{X"}$
which is a subbundle of  $N_{X"}$ : then the tangent bundle to $X$
around $X"$ is annihilated by the given subbundle $Z"$.
\item
5) $Z$ projects onto a divisor $X"\subset X$ which is focally degenerate,
$X"$ is a divisor of
$X'$ and
$X'$ is isotropically focally degenerate (again here
$\Sigma$ is not contained in ${\bf P_{\infty}}$).  Then $X$ and $X'$ are
tangent along $X"$.
\item
6) $\Sigma$ is  contained in ${\bf P_{\infty}}$, $Z$ projects
onto some affine point of $X$ , whence it dominates $X$
and $X=X'$ is obtained
via the asymptotic inverse focal construction.
\item
7) $\Sigma$ is  contained in ${\bf P_{\infty}}$, $Z$ projects
onto a component $\Delta$ of $X_{\infty}$. In this case  $Z$ is the
restriction of
$N_X$ to
$\Delta$, the second projection to ${\bf P_{\infty}}$ is not
surjective. This case is characterized by the property that
$\Delta \subset X_{\infty}$ be projectively isotropically  degenerate,
which is equivalent to the property that $\Delta$ be obtained via the
isotropic projective inverse focal construction (this case will be treated
separately in the next proposition).
\end{itemize}

Conversely, start from any  admissible pair of a variety $ \Sigma$
not contained in ${\bf P_{\infty}}$,
and of an algebraic section $r(s)$.
 Consider  the algebraic
set $X'$ obtained from the inverse construction : then
$X'$ is   focally degenerate ( if it has two components,
this means that each of them is focally degenerate)
 and isotropically focally degenerate
iff the fibres of $Z'  \rightarrow \Sigma$ are affine spaces
of dimension $m - dim \Sigma$ (then
$Z' = N_{X'}$) .

All the isotropically focally degenerate varieties $X$ are gotten
by the inverse construction as such an $X'$.

In  case 6), where $\Sigma$ is  a component of the strict
focal locus contained in ${\bf P_{\infty}}$,
all such weakly focally degenerate varieties are obtained from
the asymptotic inverse focal construction.

Let us consider the remaining cases where $\Sigma$ is
not contained in ${\bf P_{\infty}}$.

Then the weakly focally degenerate varieties in the non vertical case
(i.e., when $Z$ dominates $X$) are
gotten either\\
 (i) as a component of such an $X'$, or \\
(ii) as a divisor  in an
isotropically focally degenerate variety $X'$, which is transversal
 to the general fibres
$X'_s$ of $N_X' \rightarrow \Sigma_{X'}$
and where $dim\, \Sigma_{X'}\leq m-2$.

Instead, in the vertical case, the weakly focally degenerate varieties
are given as varieties containing a focally degenerate divisor $X"$
such that either\\
(i) $X"$ is a component of $X'$, with a component
$Z"$ of the ramification locus $Y_{X"}$ which is a subbundle of $N_{X"}$,
and such that the tangent bundle of
$X$ along $X"$ is given by the annihilator of the subbundle $Z"$ or\\
(ii) $X"$ is a divisor of $X'$, $X'$ is isotropically focally degenerate
with $\dim\, \Sigma_{X'}\leq m-2$, $X"$ is transversal to the fibres
$X'_{s}$ of $N_X' \rightarrow \Sigma_{X'}$, and $X$ and $X'$ are tangent
along $X"$.

\end{teo}

Proof\\
We discuss first of all the case where $\Sigma$ is
not contained in ${\bf P_{\infty}}$ (whence $Z$ does not
project to ${\bf P_{\infty}}$ under the first projection) .

Around each smooth point of $X$ there are a Zariski open
set $U$ of ${\bf C}^{m}$ and polynomials $ F_1(x),..    .. F_{m-n}(x)$
such that
$ X \cap U$ = $ \{ x \in U | F_1(x) = ..    .. F_{m-n}(x)=0 \}$
and such that $ X \cap U$ consists of smooth points.

Therefore, the gradients of the polynomials $ F_1(x),..    .. F_{m-n}(x)$
yield a framing of the Euclidean normal bundle on $X \cap U$,
and the endpoint map is locally given by

$\epsilon (x,\lambda_1, .. \lambda_{m-n}) = x + \Sigma_{i=1,.. m-n}
\lambda_i  \nabla F_i(x)$.

We choose as we did before a component $Z$ of the ramification locus
$Y_X$ which maps onto an irreducible variety $\Sigma$ of dimension
$ a \leq m-2$ (respectively $ a \leq m-1$ in the
focally isotropically
degenerate case) and local coordinates

$s=(s_1, ...., s_a)$ for the points of $\Sigma$ and
  $ t= (t_1, .. t_{\nu-a})$
for the fibres of $\pi$, where $\nu$ equals $m-1$ in the
non focally isotropically
degenerate case, otherwise $m = \nu$ and $Z = N_X$.

Whence, we have local functions
$ x(s,t) , \lambda(s,t) $ parametrizing the points of $Z$ ,

and a function $ O(s)$ parametrizing the image $\pi (Z) = \Sigma$
such that\\

1") $ F_i(x(s,t)) \equiv 0 $ $\forall i$\\

2") $ x(s,t) - O(s) \equiv - \Sigma_{i=1,.. m-n}
\lambda_i (s,t) \nabla F_i(x(s,t))  $.\\

Differentiating 1") with respect to both sets of variables $s,t$,
we infer that

$ < \nabla F_j( x(s,t)) ,( \partial x(s,t)/ \partial t_i )>
\equiv < \nabla F_j( x(s,t)) ,( \partial x(s,t)/ \partial s_h) > \equiv 0
 \forall i,j,h$.

By 2") $ x(s,t) - O(s) $ is a normal vector, whence\\

3") $< x(s,t) - O(s) , x(s,t) - O(s) > \equiv r(s)$.

and

4") $(\partial r(s)/ \partial s_j) = - 2 < x(s,t) - O(s), ( \partial O(s)/
 \partial s_j)>$.\\

Therefore, for fixed $s$, the projection $X_{s}$ of the fibre $Z_{s}$
(generally a manifold of dimension $\nu-a$) is contained in the
intersection $X'_{s}$ of a sphere $S_{s}$ of centre
$O(s)$ and radius $ r(s)^{1/2}$ with an affine space $\Pi_{s}$ of
codimension
$a$ (since at the general point we can assume $ \partial O(s)/
 \partial s_1), ...  \partial O(s)/
 \partial s_a)$ to be linearly independent).

Thus, the manifold $X'_s$ has dimension either
$m-1-a$ or $m-a$ ( but in the latter case, by Lemma 3,
 the orthogonal to $T\Sigma_{O(s)}$
is contained in $T\Sigma_{O(s)})$.

Consider as before the locus $X'$ given as the projection of the locus

$Z' \subset {\bf P}^{m} \times \Sigma $ defined as

3"') $ \{(x,s)| < x - O(s) , x - O(s) > \equiv r(s)$

4"') $(\partial r(s)/ \partial s_j) = - 2 < x - O(s), ( \partial
 O(s)/ \partial s_j)>\}$.

\begin{lem}
$Z'\subset N_{X'}$
\end{lem}
Proof\\
We must prove that the vector $x-O(s)$ is normal to $X'$.  This follows
>from the calculation of the tangent space to $Z'$ at the point $(x,s)$
that we have done above (cfr. 5"').
$\Box$

\begin{cor}
Each component of $X'$ is focally degenerate and indeed
 isotropically focally degenerate iff
$X'_{s}=\Pi_{s}$ (whence, in the latter case, $X'$ is also irreducible).
\end{cor}
Proof\\
 If $X'_{s}=\Pi_{s}$, then
$Z'$ is irreducible and $dim\, Z'=dim\, N_{X'}=m$ so that $Z'=N_{X'}$
and $X'$ is irreducible and isotropically focally degenerate.

If $dim\, X'_{s}=m-1-a$ (in this case $a\leq m-2$), then $Z'$ is a divisor
in $N_{X'}$ and hence
$\Sigma$ is  contained in $\Sigma_{X'}$.
Either $\Sigma$ is a component of $\Sigma_{X^0}$,
for each component $X^0$ of $X'$, or there is a
 component $X^0$ of $X'$ which  is focally
isotropically degenerate.

Assume that the latter holds:   then, for general $O(s) \in \Sigma$,
 $X'_s$  is a divisor of the fibre
of $ N_{X^0}  \rightarrow \Sigma_{X^0}$, whence by dimension
reasons $ \Sigma = \Sigma_{X^0}$.

Since the direction of $\Pi_{s}$ is the vector subspace $ W = T\Sigma_{O(s)}
^{\perp}$, and $ \Sigma = \Sigma_{X^0}$, it follows that
$ N_{X^0_s}  = \Pi_{s}$.

Moreover, being $X^0$ isotropically focally degenerate,
by lemma $3$ follows that $W$ is totally isotropic,
whence the quadratic function
$< x- O(s), x- O(s) >$ is then constant on $\Pi_{s}$, contradicting
the fact that for general $s$ $X'_{s}$ is a nonempty and proper
divisor in $\Pi_{s}$.
$\Box$\\

If $X$ is focally isotropically
degenerate, the projection $X_s$ of $Z_s$  has dimension $m-a$,
whence it equals  $X'_s$, and it follows
immediately that $X$ equals $X'$.

Suppose then that $X$ is not isotropically focally degenerate, and let
$X''$ be the projection of $Z$, that is the closure of $\cup _s X_s$. Thus
$X"\subset X$ and $X"\subset X'$.\\
  Assume first that $dim\, X'_{s}=dim\, X_{s}=m-1-a$. Therefore, $X"$ is a
component of
$X'$ and $Z$ equals  a component $Z"$ of $Z'$.
It follows that either $X" = X$ and case 2) of the theorem occurs, or $X"
\subset X$ would be a divisor and $Z$ would be the restriction to $X"$ of
the normal bundle
$N_X$, a subbundle of the normal bundle $N_{X"}$.

Whence, $X"$ is focally degenerate, with a component
$Z = Z"$ of the ramification locus which is a projective subbundle of
>$N_{X"}$, and case 4) occurs. Any variety $M$ containing $X"$ as a divisor,
and with tangent bundle annihilated by the given subbundle
 would be a weakly
focally degenerate variety with $\Sigma$ in the focal locus.

In other words, in the vertical case, the inverse focal construction
can by no means reconstruct $X$, but only the first order neighbourhood
of $X$ along $X"$.

Finally, there remains the case where $dim\, X_s = m-1-a ,\, dim X'_s =
m-a$, in which case $X'$ is isotropically focally degenerate.
 Then $Z$ is a divisor in $Z'=N_{X'}$.

Assume $X" = X'$ : since then $X' \subset X$, but $X' \neq X$ since
$X$ is not isotropically focally degenerate, we get that
$ Z = N_X |_{X'} \subset Z' = N_{X'}$, and we are again in case 4).

Thus we may consider the remaining cases where
 $X"$ is a divisor of $X'$.
Furthermore, either $X"=X$ or $X"$ is a divisor in $X$.  If $X=X"$, then
$Z$ is the restriction of $N_{X'}$ to $X$, and case 3) occurs.
If $X"$ is a divisor of $X$, we have that $Z=N_{X}|_{X"}=N_{X'}|_{X"}$
so that
$X$ and $X'$ are tangent along $X"$ and case 5)  occurs.
Conversely, let $X'$ be a isotropically focally degenerate variety and let
$X$ be a divisor inside $X'$; since $ N_{X'}|_{X} \subset N_X$ is a
divisor, it follows immediately that, setting $Z= N_{X'}|_X$, the image of
$Z$ is contained in $\Sigma_{X'}$. If moreover, as it should be, the
divisor $X$ is transversal to the fibres $X'_s$, then its image equals
$\Sigma_{X'}$, whence $Z$ will make $X$  weakly focally degenerate if and
only if $dim\,
\Sigma_{X'}
\leq m-2$. More generally, if $M$ is any variety containing $X$
as a divisor and such that $M$ and $X'$ are tangent along $X$, then
$M$ is weakly focally degenerate.\\

Let us then consider case 6) : then, analogously to the case of
hypersurfaces we can find a parametrization
 $O(s)$  of $\Sigma$ in homogeneous coordinates such that

 $ O(s) \equiv - \Sigma_{i=1,.. m-n}
\lambda_i (s,t) \nabla F_i(x(s,t))  $.\\

Then  $  < \partial x (s,t)/ \partial t_i  , O(s) > \equiv
< \partial x (s,t)/ \partial s_j  , O(s) > \equiv 0$.

>From the first equalities we conclude that there exists a local function
$r(s)$ such that

6)  $  < x (s,t) , O(s) > \equiv r(s)$.

One moment's  reflection, since the vector $O(s)$ gives homogeneous
coordinates for $\Sigma$, shows that indeed $r(s)$ globalizes
to a (possibly multivalued and with poles) section of
${\cal O}_{\Sigma}(1)$.

>From the second equalities follows also

7) $(\partial r(s)/ \partial s_j)\equiv < x(s,t) , ( \partial
 O(s)/ \partial s_j)>, for  j = 1,... a$.\\

Thus an entirely similar argument yields that $X$ is gotten
>from the asymptotic inverse focal construction, and conversely
if $X$ is obtained in this way then $X$ is weakly
focally degenerate and we are in case 6).

 $\Box$\\

Let us discuss case 7),  where the whole condition of degeneracy
 bears on $X_\infty$,
and tells that, $O(s)$ being the $V'$- vector valued function giving
local homogeneous coordinates around a smooth point
 of $\Sigma$ as usual,
there is a local function $\lambda (s,t)$ and a local
parametrization $x(s,t)$,  of $X_\infty$ this time,
and giving homogeneous coordinates, such that

$\lambda (s,t) x(s,t) - O(s) $ is a normal vector to $X_\infty$
at the point $x(s,t)$, in the sense that

$ <\lambda (s,t) x(s,t) - O(s), x(s,t)> \equiv \\
< \partial x (s,t)/ \partial t_i  ,\lambda (s,t) x(s,t) - O(s) >
 \equiv \\
< \partial x (s,t)/ \partial s_j , \lambda (s,t) x(s,t) - O(s) >
 \equiv 0 $.

At the points where $\lambda (s,t)$ is not vanishing we can
replace the parametrization $x(s,t)$ by  $\lambda (s,t) x(s,t)$,
so with these new homogeneous coordinates we have

  $ I)  < x(s,t) - O(s), x(s,t)> \equiv \\
II)  <\partial x (s,t)/ \partial t_i  , x(s,t) - O(s) >
 \equiv \\
III)  <  \partial x (s,t)/ \partial s_j ,  x(s,t) - O(s) >
 \equiv 0 $.

Deriving  equation I) with respect to $\partial/ \partial t_i $, and using
II) we obtain

$ IV) <\partial x (s,t)/ \partial t_i  , x(s,t)> \equiv 0$

whereas, applying  $\partial/ \partial s_j $ to I) and using III)
we get

$ V) <  \partial x (s,t)/ \partial s_j ,  x(s,t)> \equiv
<  \partial O (s)/ \partial s_j ,  x(s,t)>. $

IV) yields

$ A) <x(s,t), x(s,t)> \equiv <O(s), x(s,t)>\equiv r(s) $ which implies,
together with V) :

$ B) <  \partial O (s)/ \partial s_j ,  x(s,t)> \equiv 1/2
\partial r (s)/ \partial s_j $.

Since we chose a smooth point of $\Sigma$ the $a+1$ vectors

$ O(s), \partial O (s)/ \partial s_1 ,  ... \partial O (s)/ \partial s_a $
are linearly independent, and it follows that the vectors
$ x(s,t)$ , for $s$ fixed,
vary in an affine space $X"_s$ of dimension $m-1-a$.

Since however $X_s$ is assumed to have dimension exactly equal
to $ m-1-a$, it follows that $X_s = X"_s$,
where  $X"_s$ is defined by the equations

$ A') <O(s), x >  \equiv r(s) $

$ B') <  \partial O (s)/ \partial s_j ,  x> \equiv 1/2
\partial r (s)/ \partial s_j $.

However, also the equality $ < x, x> \equiv r(s) $ must be satisfied
on $X_s = X"_s$, thus by Lemma 3 we get an affine linear subspace
with direction $W$ which is totally isotropic, and is contained
in the orthogonal $W^{\perp}$  to  $W$.

The conclusion is that the projective tangent space to $\Sigma$
at any smooth point
has a totally isotropic  annihilator .

\begin{df}
 Let $\Sigma$ be a projective subvariety of the projective
space  ${\bf P}(V')={\bf P_{\infty}}$ associated to a
vector space $V'$ of dimension $m$ endowed with a non degenerate
quadratic form $Q_\infty$, such that any point $O(s)$ of
$\Sigma$ the projective tangent space to $\Sigma$ at $O(s)$
(a vector subspace of $V'$) has a totally isotropic annihilator.

Let $r(s) $ be an algebraic section of ${\cal O}_{\Sigma}(1)$
and consider the developable variety $X"$ defined by the union of
the subspaces $X"_s$  defined by the equations A') and B').

Assume moreover that $r$ be $\bf admissible$ in the sense that
the local function (constant on $X"_s$)

$<x(s,t), x(s,t)> \equiv  R(s) $ be equal to $r(s)$.

Then we shall say that $X"$ is projectively isotropically
 degenerate and that $X"$ is obtained via the isotropic
projective inverse  focal construction from the
admissible pair $(\Sigma, r)$.

\end{df}

\begin{prop}
Assume that $X$ is weakly focally degenerate and that a component
 $\Sigma$  of the focal locus is  contained in ${\bf P_{\infty}}$,
with the corresponding component $Z$ of $Y_X$ projecting
onto a component $\Delta$ of $X_{\infty}$ ( case 7) of theorem 3).
 In this case  $Z$ is the restriction of $N_X$ to
$\Delta$,  $\Delta$
is projectively isotropically  degenerate. Conversely, if  $\Delta$
is obtained via the isotropic
projective inverse focal construction , then $X$ is weakly focally
degenerate and we are in case 7) of theorem 3.

\end{prop}
Proof\\
If $X$ is as in case 7) of theorem 3, then we have already seen that
$\Delta \subset X_{\infty}$ is projectively isotropically  degenerate.

It remains to prove the converse, which follows since A'), B') and
our assumption $R(s) = r(s)$ imply A), B) by which immediately
follow I), II), and III), whence $x(s,t) - O(s)$ is a normal vector
to $X_{\infty}$. Since $X"=\Delta$ and $X"_s$ has dimension $m-1-a$
we get a component $Z$ of dimension $m-1$ projecting onto the
$a$-dimensional variety $\Sigma$ contained in the hyperplane
at infinity and we are done.

$\Box$\\

\begin{rem}
It follows from the previous theorem that any variety $\Sigma$ is
a component of some focal locus.

Moreover, in the asymptotic inverse construction, we see
immediately that the tangent space at a point of $X_s$ depends
only upon $s$ , so that then  our $X$ is developable.

In particular, if $m=3$ we get either a linear subspace or a developable,
whence singular, surface.

Observe finally that if $X$ is orthogonally general and projective,
then only cases 6) or 7) can a priori occur.

For case 6),
start choosing $X_\infty$ as a smooth and transversal
variety to $Q_\infty$,
 apply then the asymptotic inverse focal construction : then we get
a variety $X$ which will be orthogonally general exactly iff
$X$ is smooth. But the smoothness of $X$, as we have just seen,
is the main obstruction.
\end{rem}
\begin{ex}
Let $m=3$, and let $\Sigma$ be the line at infinity parametrized as

$ O(s) = (0,0,1,s)$, and set, in these affine coordinates,
$r(s) = s^2/2$.

Then an easy computation for the asymptotic inverse focal
construction yields the quadric cone

$ X = \{ x | x_0 x_2 - x^2_3/2 = 0 \}$, whose vertex lies at
infinity.

If instead we choose $ O(s) = (0,1,s,s^2)$, and
$r(s) \equiv 1$, $X$ will be the quadric cone

$ X = \{ x | 4 ( x_1 - x_0) x_3 - x^2_2= 0 \}$, whose vertex does
not lie at infinity.

\end{ex}

\begin{ex}
Let us now consider the most classical example, namely the
rotational torus $X$ obtained rotating a circle of  radius,
say, $1$ around the point with coordinates $(2,0)$.
This is the example of a strongly focally degenerate variety.
\end{ex}

The equation $F$ of $X$, in affine coordinates
$(x,y,z)$ for which $Q_\infty$ yields the standard Euclidean
scalar product, is then given, setting

$q(x,y,z) = ( x^2 + y^2 + z^2 + 3)$, or , in homogeneous
coordinates $(x,y,z,w)$, $q(x,y,z,w) = ( x^2 + y^2 + z^2 +
 3w^2)$, by

(*) $q^2 - 16 (x^2 +y^2)w^2$.

The intersection with the plane at infinity is precisely
our conic $Q_\infty = \{ q=w=0 \}$, which is a double
 curve for the quartic surface $X$. Moreover, $Sing(X)$
consists of $Q_\infty$ and of the two points $\{P,P'\}
 = \{ q=x=y=0\}= \{ (z^2 + 3w^2)=x=y=0\}$.

Now, a classical and easy formula for a rotation surface
of a curve
$C$ =$ r(s), z(s)$ parametrized by arclength,

$ x (s,\theta) = r(s) cos( \theta)$

$ y (s,\theta) = r(s) sin(\theta)$

$ z(s,\theta) = z(s)$

is that the two principal curvatures equal $k(s)$ , $z'(s)/r(s)$.

In this case, $ r(s), z(s) = (2 + cos (s), sin(s))$, whence
 $ k \equiv 1$ and $z'(s)/r(s)$ = $ 1 - 2/r(s)$.

These formulae are easily rationalized on our surface $X$ since
 $q^2 = 16 (x^2 +y^2)$, whence $r = q/4$.
Therefore the critical points are obtained by taking the
multiples of the unit normal by the opposites of
 their inverses, i.e., $- 1$ and $- q / (q-8)$.
Finally, the unit normal is obtained by the gradient of $F$

$ \nabla F= ( 4x(q-8) , 4y(q-8), 4qz)$ upon dividing by its norm,
which equals

$ | \nabla F| = 4 ( (q-8)^2 (x^2 + y^2) + q^2 z^2)^{1/2}$ =

$ ( 16(q-8)^2 (x^2 + y^2) + 16 q^2 z^2)^{1/2}$=
$ q ( (q-8)^2  + 16  z^2)^{1/2}$. But since

$ z^2 = q - 3 - q^2/16$ , we get  $ q (64 -48)^{1/2} = 4 q$,

and the focal locus is obtained for the values $ \lambda =
- 1/4q , \lambda = - 1/4(q-8)$ as the image of the endpoint map
$ (x,y,z) + \lambda \nabla F(x,y,z)$.

For $\lambda = - 1/4(q-8)$ we get the points
$(0,0,z (8/q-8))$ , for $ \lambda = - 1/4q $ we get the points
$(8x/q, 8y/q, 0)$.

The conclusion is that the focal locus consists of the
$z$-axis and of the circle $ z=0, x^2 + y^2 =4$. That is,
our surface is strongly focally degenerate, and we can indeed
see geometrically the two families of circles corresponding
to the two components of the focal locus.

We end this protracted example by observing that the rotation
surface is clearly a rational surface. Indeed, we can say more,
since a smooth model is obtained by blowing up the singular
conic $Q_{\infty}$ and the two points $P,P'$.

Let $R$ and $E, E'$ be the respective exceptional divisors
in the blow-up $ \widetilde{\bf P}$ of ${\bf P}^{3}$: the first
is a ruled surface ${\bf P}( {\cal O}_{\bf P^1} \bigoplus
{\cal O}_{\bf P^1} (2))$ , the other are two ${\bf P}^{2}$'s.

Let $S$ be the strict transform of $X$: it belongs to
the linear system $| 4H - 2 R - 2E - 2E'|$, whereas the canonical
system of $ \widetilde{\bf P}$ equals $| - 4H + R + 2E + 2E'|$.
Thus $S$ belongs to $ | -K - R|$ , and by the exact sequence

$ 0 \rightarrow {\cal O}_{ \widetilde{\bf P}}$
$(-S) \rightarrow {\cal O}_{ \widetilde{\bf P}}$ $ \rightarrow
{\cal O}_S \rightarrow 0$

we infer $ h^i({\cal O}_S ) = h^{i+1}
 ({\cal O}_{\widetilde{\bf P}}$ $ (K+R) = h^{2-i}
({\cal O}_{\widetilde{\bf P}}$ $(-R) = 0$,
since $R$ is irreducible.
$S$ is clearly then rational, and the anticanonical
effective divisor has self-intersection $4$.

\begin{ex}
More generally, for a rotation surface
$ (r(s) cos( \theta), r(s) sin(\theta), z(s))$

the unit normal is given by
$ (z'(s) cos( \theta), z'(s) sin(\theta), r'(s))$ ,
therefore we see easily that the focal locus consists
of the $z$-axis and of the rotation surface obtained
by rotating the evolute of the plane curve $C =\{ r(s),z(s)\}$
we were starting with.

Therefore, general rotation surfaces provide examples
of weakly but not strongly focally degenerate varieties.

\end{ex}

\begin{ex}
This last example shows the important role of the algebraic
function r(s).

Let $\Sigma$ be the line $\{ (0,0,s) \in {\bf C}^{3} \}$ : then
if we take the function $r(s) \equiv R$, where $R \in {\bf C}$
is a constant, the inverse construction yields a cylinder $X'$.
Instead, if we choose $r(s) \equiv R + s^2$, we obtain as $X'$
simply a circle in the plane $z=0$.
\end{ex}

\section{Isotropically focally degenerate hypersurfaces}

In the preceding section we gave a characterization, in terms
of the inverse focal construction, of the focally isotropically
degenerate varieties. However, in general such a construction
yields a hypersurface, which is only weakly degenerate, and
although in the next section we shall write down conditions which
characterize the focally isotropically
degenerate case,
in the case of hypersurfaces, we can give an easier characterization
for the isotropically degenerate case with a direct proof.

Let thus $F(x_1,... x_{n+1}) = 0$ be the polynomial equation
of an affine hypersurface $X$, which we may assume, without loss
of generality, to be irreducible.

Again the gradient $\nabla F$ of $F$ gives a map
of the Normal Bundle $N_X$, $\pi:N_{X}
 \rightarrow {\bf P^{n+1}}$
which we will also call the endpoint map \\

$\epsilon (x,\lambda) = x + \lambda \nabla F(x)$,
where $x$ is a point of $X$ ( thus, for $\lambda =0$
we reobtain the points of $X$).

\begin{prop}
Let $X$ be a projective hypersurface :
then $X$ is focally isotropically
degenerate if and only if $X$ coincides with its focal
locus $\Sigma_X$.
\end{prop}
Proof\\
In this case the focal locus equals the image $\Sigma_X$ of the
map $\pi : N_X \rightarrow {\bf P}^{n+1}$, and since  $X$ may
be assumed to be irreducible, $N_X$ is irreducible, whence
$\Sigma_X$ is also irreducible. But $X$ is contained in
$\Sigma_X$ and has not lesser dimension, thus equality holds.
$\Box$\\

\begin{rem}
We derive thus the equality

$ F (x + \lambda \nabla F(x) ) \equiv 0$ $ \forall \lambda$.

In particular $(d /d \lambda ) F (x + \lambda \nabla F(x) )
 \equiv 0$, and, for  $\lambda =0$, we get

(I) $ < \nabla F(x) , \nabla F(x) > \equiv 0$.
\end{rem}

By the previous proposition the general fibre of $\pi$
has dimension $1$, and for each $x_0 \in X$, $\lambda_0 \in
{\bf C}$ there exists a curve

(II) $ x(t),\lambda (t)$ such that  $ x(0)= x_0,
\lambda (0) = \lambda_0$, which is a fibre of $\pi$.

Since a fibre intersects a normal line $x_0 \times {\bf C}$
in at most one point, it follows that up to a birational
transformation we can take $ (x_0,t)$ as coordinates
on $N_X$ by taking the curves $ x(x_0,t),\lambda (x_0,t)$
satisfying (II) for $\lambda_0 = 0$,
 and assume that the curve $x(x_0,t)$ is a non constant
curve in $X$ satisfying

(III) $ x(x_0,t) + \lambda (x_0,t)\nabla F(x(x_0,t)) \equiv x_0$.

We argue as in the preceding section :

$ x(x_0,t) - x_0 \equiv - \lambda (x_0,t)\nabla F(x(x_0,t))  $,

thus by (I) our usual function $r(x_0) \equiv 0$ and

(IV) $ < x(x_0,t) - x_0 , x(x_0,t) - x_0  > \equiv 0$.

In this case we also get, if $s=(s_1, ...., s_n)$ are local
coordinates for $x_0 \in X$, that  $(d r(s)/ d s_j)\equiv 0 $ and

(V)  $ 0  = - 2 < x(x_0 (s),t) - x_0(s),
( d x_0(s)/ d s_j)>$ for each $s,t$.

Since the tangent space to $X$ at $x_0$ has dimension $n$,
we infer that, fixing $s$ and varying $t$, we obtain a curve $x(x_0,t)$
which moves on the line
through $x_0$ with direction $\nabla F(x_0)$.

We can thus write

(VI) $x(x_0,t) = x_0 + \mu (x_0,t)\nabla F(x_0) $,

and then (VI) and (III) combine to yield

(VII) $ \lambda (x_0,t) \nabla F(x_0,t) \equiv
- \mu (x_0,t)\nabla F(x_0) $.

Since the function
$\lambda (x_0,t)$ is non zero, it follows that
not only the line
through $x_0$ with direction $\nabla F(x_0)$ is contained
in $X$, but also that the normal direction stays constantly
proportional to $\nabla F(x_0)$ on it.

We have thus proven the following

\begin{teo}
A hypersurface $X$ is isotropically focally degenerate if
and only if it is isotropically developable, i.e.,
for each point the normal line is contained in $X$, and
along this line the tangent space to $X$ does not vary.
\end{teo}

We would like now to give some examples and show where
lies the difficulty in the fine classification of
isotropically focally degenerate hypersurfaces.

It is classically known that in 3-space the
analytical surfaces which are developable are only
cones, cylinders, and tangential developable surfaces.

\begin{prop}
Assume $X$ is an isotropically developable surface.
Then ,  if $X$ a cylinder then $X$ is a plane.
If $X$ is a cone , it is the cone over $Q_{\infty}$
with vertex in a point of affine space.
\end{prop}
Proof\\
If $X$ is a cylinder, then the generatrices are the normal
lines, therefore the normal direction is constant on
the whole surface and the surface is a plane.
If $X$ is a cone, with vertex, say, at the origin,
then the vectors $x$ and $\nabla F(x)$ are proportional,

but the vector $\nabla F(x)$ is always isotropic, whence
$ < x, x> \equiv 0$ on $X$, q.e.d.
$\Box$\\

Let us now discuss the tangential surface $X$ of a curve $C$.

We write as usual $X$ parametrically as

$ x(s,t) = \alpha (s) + t \alpha '(s)$,

so that the tangent plane is generated by the two vectors

$  \alpha '(s), \alpha ''(s)$.

Up to local analytic reparametrization we can assume
that one and only one
of the following two possibilities occurs:

(I) $ < \alpha '(s), \alpha '(s)> \equiv 0$\\
(U) $ < \alpha '(s), \alpha '(s)> \equiv 1$.

In both cases follows that

(*) $ < \alpha '(s), \alpha ''(s)> \equiv 0$.

In the isotropic case (I), then clearly $\alpha '(s)$
is a normal vector to $X$, constant on the generatrices,
and our $X$ is thus isotropically developable.

We could stop our discussion here, since the isotropic
ruling, in the situation we are interested in, is obtained
by fixing $s$ and varying $t$, which means that we are
in principle through with our discussion.
Nevertheless, for curiosity, we analyze also the unitary case
which we could avoid to consider in view of the assumption
that our surface is not only developable, but also
isotropically developable.

\begin{lem}
The unitary case (U)  occurs only if the curve
$C$ is a plane curve, thus its tangential surface
is a plane.
\end{lem}

Proof\\
In the unitary case (U), the normal vector must be proportional
to $\alpha ''(s)$, whence $X$ is isotropically
degenerate if and only if

(**) $ < \alpha ''(s), \alpha ''(s)> \equiv 0$.
Now, by taking
 derivatives of (*) and (**), and using (**), we obtain

(***) $ < \alpha ''(s), \alpha '''(s)> \equiv 0$\\
		$ < \alpha '(s), \alpha '''(s)> \equiv 0$

>from which it follows that $ \alpha'''(s), \alpha ''(s)$
are proportional vectors, whence also

$ < \alpha '''(s), \alpha '''(s)> \equiv 0$.

By induction, we show that for each integer n

(*n*) $ < \alpha ''(s), \alpha ^{(n)}(s)> \equiv 0$\\
		$ < \alpha '(s), \alpha ^{(n)}(s)> \equiv 0$\\

whence $  \alpha ''(s), \alpha ^{(n)}(s)$ are proportional
and thus also

$ < \alpha ^{(n)}(s), \alpha ^{(n)}(s)> \equiv 0$.

Consider now the Taylor development of $\alpha (s)$
at any point : from the fact that all higher derivative
vectors are proportional follows that $\alpha (s)$ yields
a plane curve.

But this means that its tangential surface is a plane.

$\Box$\\

It is now clear that in order to classify the non-trivial
isotropically developable surfaces in 3-space we would need
to classify the isotropic space curves $C$ ( i.e.,
those whose tangent vector
is always an isotropic vector,that is, (I) holds).

Now, the condition that $C$ is algebraic is an obstacle!

Indeed, $C$ will be the  birational image of a smooth
curve $B$ , given through $4$ sections  $(s_0,s_1,s_2,s_3)$
of a line bundle on $B$: the isotropicity condition
amounts to the following equation ( where $'$ represents
the derivative with respect to a local parameter on $B$)

($\cal E$)
$ \Sigma_{i=1,2,3} ( s'_i s_0 - s'_0 s_i )^2 \equiv 0$.

\begin{ex}
It is rather easy to give examples of rational curves which
are isotropic.

It suffices, chosen an affine coordinate $t$ on $\bf C$,
to set $s_0 \equiv 1$ and let $(s_1,s_2,s_3)$ be polynomials
in $t$ such that their derivatives satisfy

$ \Sigma_{i=1,2,3} ( s'_i )^2 \equiv 0$.

In other words ,$(s'_1,s'_2,s'_3)$ give a rational
parametrization of the conic $Q_{\infty}$ and
$(s_1,s_2,s_3)$ are taken to be the integrals of
the three polynomials $(s'_1,s'_2,s'_3)$.

In this way we see more generally that, up to translation,
our curve $C$ is determined by our map $B \rightarrow
Q_{\infty}$.

Using in our particular case of the rational curves
 fixed isomorphisms of $ \bf P^1$ with $B$ and with
$Q_{\infty}$, we obtain that our isotropic rational
curves are parametrized by a
pair of  polynomials $f_0(t),f_1(t)$.

In concrete terms , we may take

$(s'_1,s'_2,s'_3)$ = $(f^2_0 + f^2_1, if^2_1 -i f^2_0,
 -2i f_0 f_1)$.

Assume now not only that the map $f_0(t),f_1(t)$ is of positive
degree and is primitive (does not factor through an
intermediate cover), e.g. it could be a cyclic Galois cover
of prime order $p$.


If the map $(s_1,s_2,s_3)$ would not be birational onto
its image, then the tangent map from $C$ to $Q_{\infty}$
would be a birational isomorphism.

But, in the example we gave above, $f_0(t) = 1,f_1(t) = t^p$,
we see immediately that $(s_1,s_2,s_3)$ are not polynomials
in $ t^p$.


\end{ex}

\section{Isotropically focally degenerate varieties and further examples }

In the previous section we have given a classification, and concrete
examples of isotropically focally degenerate hypersurfaces.

It is easy to obtain concrete examples in higher codimension by
the following simple device : consider varieties $M \subset  \bf C^m$,
$W \subset  \bf C^w$ and consider the product variety $ X = M \times W$
in the orthogonal direct sum $ \bf C^m \bigoplus ^{\perp} \bf C^w =
 \bf C^{m+w}$ .

It is immediate to see that in this case the normal bundle of $X$
is a product,  likewise the endpoint map.

\begin{rem}
If  thus $M \subset  \bf C^m$ and $W \subset  \bf C^w$ are
isotropically focally degenerate, then $ X = M \times W \subset
 \bf C^{m+w}$ is also isotropically focally degenerate, and
$\Sigma_X = \Sigma_M \times \Sigma_W $. In particular, we obtain
in this way $\Sigma_X$ of arbitrary codimension.

We obtain also, by letting $M$ be an  isotropically developable hypersurface,
and $W$ general,
an example of a variety $X$ of arbitrary codimension  which
is isotropically focally degenerate, and whose $\Sigma_X$ is a hypersurface.
\end{rem}

We now finally observe that the inverse focal construction
 gives a characterization
of the isotropically focally degenerate varieties in terms of their
focal variety $\Sigma_X$.

\begin{teo}
Let $\Sigma$ be a projective variety of dimension $a$, and let $r(s)$
be an algebraic function on its affine part.
Assume moreover that

1) at any point $O(s)$ of
$\Sigma$ the vector tangent space to $\Sigma$ at $O(s)$
(a vector subspace of $V'$) has a totally isotropic annihilator.

Then, if $X$ is gotten from $(\Sigma, r)$ via the inverse focal construction,
and moreover

2) the algebraic function $r(s)$ satisfies the conditions

2.1) $dr(s) \in Im ( T_{\Sigma, s} \rightarrow ^{Q_\infty}\rightarrow
T_{\Sigma, s}^{\vee })$

2.2) given $\xi$ with Im($\xi$) = $df$ , then $1/4 < \xi, \xi > = r(s)$

then $X$ is isotropically focally degenerate and $\Sigma = \Sigma_X$.

\end{teo}

Proof\\
This follows immediately from Lemma 3, since conditions 1)
and 2.1) imply that on the affine space given by
equations 4"") the quadratic function $Q_\infty$ is constant,
and 2.2) guarantees that this constant equals $r(s)$, whence also 3"")
is satisfied and thus the sphere $X_s$,
fibre over  the point $O(s)$, is then an affine space of dimension
$m-a$.

$\Box$ \\

\begin{rem}
The above theorem immediately implies the characterization given
in the previous section of the isotropically focally degenerate
 hypersurfaces. Because in the case of hypersurfaces we noticed
that $X = \Sigma_X$, and then the tangential condition on $\Sigma_X$
reads out as the condition that the normal vector is isotropic,
moreover by the inverse focal construction $X$ is developable,
and the fibre dimension equals $m-a = m -(m-1)= 1$.
Thus $X$ is developable with the ruling by lines given by the
normal direction.

\end{rem}

We end by showing an explicit example of the situation considered
in case 3) of Theorem 3.

\begin{ex}
Consider first $X' \subset \bf C^6$ obtained as the product
$ X' = M \times W $ of
>two isotropically developable surfaces :

thus $X'$ has a parametrization

$( \alpha(s) + t \alpha'(s), \beta(\sigma) + \tau \beta'(\sigma)) $.

Inside $X'$ we consider the divisor $X$ obtained by setting $ \tau = t$.

Whence $X$ has a parametrization

$( \alpha(s) + t \alpha'(s), \beta(\sigma) + t \beta'(\sigma)) $,

and, remembering that $\alpha'(s) , \beta'(\sigma) $ are isotropic vectors
it follows that the normal space to $X$ is spanned, at the smooth points
of $X$ , by the three vectors $( \alpha'(s) , 0)$ , $( 0 , \beta'(\sigma) )$
and

$(- <\beta'''(\sigma), \beta'(\sigma)>  [ \alpha''(s) + t \alpha'''(s)] ,
 < \alpha''' (s), \alpha'(s)  > [\beta''(\sigma) + t \beta'''(\sigma)] )$.

The endpoint map is given by

$( \alpha(s) + t \alpha'(s) + \lambda_1 \alpha'(s) - \lambda_3
<\beta'''(\sigma), \beta'(\sigma)>  [ \alpha''(s) + t \alpha'''(s)],
\beta(\sigma) +  t \beta'(\sigma) + \lambda_2 \beta'(\sigma) +
\lambda_3  < \alpha''' (s), \alpha'(s)  >[\beta''(\sigma) + t
\beta'''(\sigma)] ) $

thus its image equals the image of the map

$( \alpha(s)  + \lambda_1 \alpha'(s) - \lambda_3
<\beta'''(\sigma), \beta'(\sigma)>  [ \alpha''(s) + t \alpha'''(s)],
\beta(\sigma)  + \lambda_2 \beta'(\sigma) +
\lambda_3  < \alpha''' (s), \alpha'(s)  >[\beta''(\sigma) +t
 \beta'''(\sigma)] ) $.

To simplify the discussion we may assume $<  \alpha''(s) ,  \alpha''(s) >
\equiv
-1 $, and similarly $ <\beta''(\sigma), \beta''(\sigma)> \equiv -1$,
therefore our formula simplifies to

$( \alpha(s)  + \lambda_1 \alpha'(s) - \lambda_3
 [ \alpha''(s) + t \alpha'''(s)],
\beta(\sigma)  + \lambda_2 \beta'(\sigma) +
\lambda_3  [\beta''(\sigma) +t
 \beta'''(\sigma)] ) $

and  we see that the image of the normal
bundle $ N_X$ will in general be dominant.

Therefore $X$ is weakly focally degenerate, not focally isotropically
degenerate, but the inverse focal construction reconstructs only
the isotropically focally degenerate fourfold $X'$.
\end{ex}

{\bf References}\\
$[Arnold et al.]$ V. Arnold, A. Varchenko, S. Goussein-Zade',
{\em Singularites des Applications Differentiables. I : Classification
des points critiques, des caustiques et des fronts d' onde},
Editions MIR, Moscow ( 1986), translation of the Russian 1982 edition.\\
$[Coolidge]$ J.L. Coolidge, {\em A Treatise on Algebraic Plane Curves},
Oxford University press (1931),
reedited by Dover, New York ( 1959).\\
$[D-F-N]$ B.A. Dubrovin, A.T. Fomenko, S.P. Novikov
{\em Modern Geometry - Methods and Applications. Part 2:
The geometry and topology of manifolds},
reedited in Springer GTM 104, (1985).\\
$[Fantechi]$ B. Fantechi, {\em The Evolute of a Plane Algebraic Curve},
UTM 408, University of Trento, (1992).\\
$[Hartshorne]$ R. Hartshorne, {\em Algebraic Geometry}, Springer
GTM 52  (1977).\\
$[Milnor]$ J. Milnor, {\em Morse Theory}, Princeton, Princeton Univ. Press,
(1973). \\
$[Salmon-Fiedler]$ G. Salmon, W. Fiedler, {\em Analytische Geometrie des
Raumes}, 2 vols., Leipzig, Teubner-Verlag, (1879-1880).\\
$[Trifogli]$ C. Trifogli, {\em Focal Loci of Algebraic Hypersurfaces: a
General Theory}, Geom. Dedicata 70: 1-26 (1998).\\
$[Trif2]$ C. Trifogli, {\em The Geometry of Focal Loci},
Ph. D. Thesis, Universita' di Milano (in preparation).\\

Address of the authors:

Fabrizio Catanese,

Mathematisches Institut der
Georg August Universit\"at G\"ottingen,\\
Bunsenstrasse 3-5, D 37073 G\"ottingen,

catanese@uni-math.gwdg.de\\

Cecilia Trifogli,

Dipartimento di Matematica F. Enriques, Universita' di Milano,\\
via C. Saldini 50,I  20133 Milano

trifogli@vmimat.mat.unimi.it
and : All Souls College, Oxford OX1 4AL

cecilia.trifogli@all-souls.ox.ac.uk
\end{document}